\documentclass[11pt]{article}
\usepackage{amsmath}
\usepackage{amssymb,amsmath}
\DeclareMathOperator{\intd}{d}

\let\pa=\partial

\let\e=\varepsilon

\let\f=\frac
\let\p=\psi
\let\om=\omega

\let\Om=\Omega

\def\cC{{\cal C}}

\def\cF{{\cal F}}

\def\cP{{\cal P}}

\def\cS{{\cal S}}

\def\cZ{{\cal Z}}

\def\no{\noindent}
\def\na{\nabla}

\def\p{\partial}

\def\dv{\mbox{div}}

\def\eqdefa{\buildrel\hbox{\footnotesize def}\over =}

\def\C{\mathop{\bf C\kern 0pt}\nolimits}
\def\DD{\mathop{\bf D\kern 0pt}\nolimits}
\def\K{\mathop{\bf K\kern 0pt}\nolimits}
\def\N{\mathop{\bf N\kern 0pt}\nolimits}
\def\Q{\mathop{\bf Q\kern 0pt}\nolimits}
\def\R{\mathop{\bf R\kern 0pt}\nolimits}

\newcommand{\Z}{{\mathbf Z}}

\newcommand{\ef}{ \hfill $ \blacksquare $ \vskip 3mm}

\newcommand{\beq}{\begin{equation}}
\newcommand{\eeq}{\end{equation}}
\newcommand{\ben}{\begin{eqnarray}}
\newcommand{\een}{\end{eqnarray}}
\newcommand{\beno}{\begin{eqnarray*}}
\newcommand{\eeno}{\end{eqnarray*}}

\renewcommand{\theequation}{\thesection.\arabic{equation}}

\newtheorem{Theorem}{Theorem}[section]
\newtheorem{Definition}[Theorem]{Definition}
\newtheorem{Proposition}[Theorem]{Proposition}
\newtheorem{Lemma}[Theorem]{Lemma}

\newtheorem{Remark}[Theorem]{Remark}

\begin{document}
\title{Global well-posedness of compressible Navier-Stokes equations for some classes of large initial data}
\author{Chao Wang $^\dag$,  Wei Wang $^{\ddag}$ and Zhifei Zhang
$^{\ddag}$\\[2mm]
{\small $ ^\dag$ Academy of Mathematics $\&$ Systems Science, CAS,
Beijing 100190, China}\\
{\small E-mail: wangchao@amss.ac.cn}\\[2mm]
{\small $ ^\ddag$ LMAM and School of  Mathematical Science, Peking University, Beijing 100871, China}\\
{\small E-mail: wangw07@pku.edu.cn,\, zfzhang@math.pku.edu.cn}}

\date{November 13, 2011}
\maketitle

\begin{abstract}
We prove the global well-posedness of three dimensional compressible Navier-Stokes
equations for some classes of large initial data, which is of large
oscillation for the density and large energy for the velocity.
The structure of the system (especially, the effective viscous flux) is fully exploited.
\end{abstract}

\renewcommand{\theequation}{\thesection.\arabic{equation}}
\setcounter{equation}{0}
\section{Introduction}

We consider the three dimensional compressible Navier-Stokes equations
\begin{equation}\label{eq:cNS}
\left\{
\begin{array}{ll}
\p_t\rho+\textrm{div}(\rho u)=0,\\
\p_t(\rho u)+\textrm{div}(\rho u\otimes u)-\mu\Delta u-\lambda\na \textrm{div}u+\na P(\rho)=0, \\
(\rho,u)|_{t=0}=(\rho_0,u_0).
\end{array}
\right.
\end{equation}
Here $\rho(t,x)$ and $u(t,x)$ denote the density and velocity of the
fluid respectively. The Lam\'{e} coefficients $\mu$ and $\lambda$ satisfy
\ben\label{assu:coeff}
\mu>0\quad \textrm{and} \quad \lambda+\mu>0,
\een
which ensure that the operator $\mu\Delta u+\lambda\na \textrm{div}u$ is elliptic.
The pressure $P$ is a smooth function of $\rho$ and satisfies
\ben
&P(0)=0,\quad P'(\overline{\rho})>0,\\
&(\rho-\overline{\rho})(P(\rho)-P(\overline{\rho}))>0
\quad \textrm{for}\quad \rho\neq \overline{\rho},\,\rho\in [0,2c_0^{-1}],\label{ass:pressure}
\een
where $0<{c_0}<\overline{\rho}<c_0^{-1}$ are fixed numbers.

The local existence and uniqueness of smooth solution of the system (\ref{eq:cNS})
was proved by Nash \cite{Nash} for smooth initial data without vacuum.
In a seminal paper \cite{Mat}, Matsumura and Nishida proved the global existence and uniqueness of smooth solution
for smooth initial data close to equilibrium. In general, whether smooth solution
blows up in finite time is an open problem. However, Xin \cite{Xin} proved that
smooth solution will blow up in finite time if the initial density has compact support.
Recently, Sun, Wang and Zhang \cite{Sun} showed that smooth solution does not blow up
if the concentration of the density does not occur.

The global existence of weak solution was proved by Hoff \cite{Hoff-JDE, Hoff-JMFM} for the discontinuous initial data with small energy.
For the large initial data, the global existence of weak solution was proved by Lions \cite{Lions}
for the isentropic Navier-Stokes equation, i.e. $P=A\rho^\gamma$ for $\gamma\ge \f95$. Feireisl, Novotn\'{y} and
Petzeltov\'{a} \cite{Fei} improved Lions's result to $\gamma>\f32$.
For the spherically symmetric or axisymmetric initial data,
Jiang and Zhang \cite{Jiang1,Jiang2} proved the global existence of weak solution for any $\gamma>1$.
However, the question of the regularity and uniqueness of
weak solution is completely open even in the case of two dimensional space.

Motivated by Fujita and Kato's result for the incompressible Navier-Stokes equations \cite{Fuj-Kat},
Danchin in a series of papers \cite{Dan-Inve,Dan-CPDE, Dan-ARMA, Dan-CPDE07} proved the local well-posedness of (\ref{eq:cNS})
for the initial data $(\rho_0-\overline{\rho}, u_0)$ belonging to the critical Besov spaces $\dot B^{\f 3p}_{p,1}\times \dot B^{\f3p-1}_{p,1}$ for $p<6$,
and global existence for the initial data close to  equilibrium in $\dot B^{\f 32}_{2,1}\times \dot B^{\f12}_{2,1}$.
This seems the largest space in which the system (\ref{eq:cNS}) is well-posed.
Indeed, Chen, Miao and Zhang \cite{CMZ-Arxiv} proved the ill-posedness of (\ref{eq:cNS})
in $\dot B^{\f 3p}_{p,1}\times \dot B^{\f3p-1}_{p,1}$ for $p>6$.

Cannone, Meyer and Planchon \cite{Can, CMP} proved the global existence of the solution
for the incompressible Navier-Stokes equations in the Besov space of negative regularity index
$\dot B^{\f 3p-1}_{p,\infty}$ for $p>3$, see also \cite{Bah} for a new proof based on the smoothing effect of the heat equation.
A very interesting point of this result is that
it allows to construct global solution for the highly oscillating initial velocity
that may have a large norm in $\dot H^\f 12$ or $L^3$. A typical example is
\begin{align*}
u_0(x)=\sin\bigl(\f {x_3} {\varepsilon}\bigr)(-\p_2\phi(x), \p_1\phi(x),0)
\end{align*}
where $\phi\in \cS(\R^3)$ and $\varepsilon>0$ is small enough. Later,
Chemin, Gallagher, and Paicu \cite{Chemin-A,Chemin-Ann} find more classes of initial data with large energy
generating the global solution. And in the proof of \cite{Chemin-Ann}, the special structure of the equation is used. We refer to \cite{Chemin-JR,Chemin-TAMS, Chemin-CMP, Paicu, Zhang-CMP}
and references therein for more relevant results.

It is highly non trivial to generalize Cannone-Meyer-Planchon's result to (\ref{eq:cNS}),
since it is a hyperbolic-parabolic coupled system. Very recently, important progress has been made by
Chen, Miao and Zhang \cite{CMZ-CPAM} and Charve and Danchin \cite{Charve},
where they construct the global solution of (\ref{eq:cNS})  for the highly oscillating initial veloicty
by proving the global well-posedness of the system (\ref{eq:cNS}) in the critical Besov space with $p>3$.

However, for all the above global existence results, the initial density is required to be
close to a positive constant in $L^\infty$ norm, which precludes the oscillation of the density at any
point. Recently, Zhang \cite{Zhang} proved the global existence and uniqueness of (\ref{eq:cNS})
for the initial density $\rho_0$ close to a positive constant in $L^2$ norm and $u_0$ small in $L^p$ norm for $p>3$,
which allows the density to have large oscillation on a set of small measure. Similar result
has also been obtained by Huang, Li and Xin \cite{Huang} for the initial data with vacuum,
but a compatibility condition is imposed on the initial data.

The proof of \cite{Zhang} and \cite{Huang} is based on Hoff's energy method, where the structure of the system
is fully exploited, especially the effective viscous flux plays an important role.
While, the proof of \cite{CMZ-CPAM}  and \cite{Charve} is based on the analysis for the linearized system, where
Harmonic analysis tools (especially Littlewood-Paley theory) are exploited,
but the structure of the system is not exploited.

In this paper, we combine two methods to find a new class of large initial data generating the global solution
of the system (\ref{eq:cNS}). This class of initial data may have large
oscillation for the density and large energy for the velocity. More precisely, we prove

\begin{Theorem}\label{thm:global}
Assume that the initial data $(\rho_0,u_0)$ satisfies
\ben \rho_0-\overline{\rho}\in H^s\cap
\dot B^{\f 3 p}_{p,1},\quad c_0\le \rho_0\le c_0^{-1}, \quad u_0\in H^{s-1}\cap
\dot B^{\f 3p-1}_{p,1} \een for some $p\in (3,6)$ and $s\ge 3$.
There exist two constants $c_1=c_1(\lambda,\mu, c_0,\overline{\rho})$ and
$$c_2=\f 1 {C\big(1+\|\rho_0-\overline{\rho}\|_{\dot B^{\f 3p}_{p,1}\cap H^s}+\|u_0\|_2\big)^{5}},\quad C=C(\lambda,\mu, c_0,\overline{\rho})$$
such that if
\ben
\|\rho_0-\overline{\rho}\|_{L^2}\le c_1,\quad
\|u_0\|_{\dot H^{-\delta}}+\|u_0\|_{\dot B^{\f3p-1}_{p,1}}\le c_2,
\een
for some $\delta\in (1-\f 3p,\f3p)$, then there exits a unique global solution $(\rho, u)$ to the compressible Navier-Stokes system
(\ref{eq:cNS}) satisfying
\beno
&\rho\ge \f {c_0} 4,\quad \rho-\overline{\rho}\in
C([0,+\infty);H^s),\\
&u\in C([0,+\infty);H^{s-1})\cap L^2(0,T;H^{s})\quad\textrm{ for any }\quad T>0. \eeno
\end{Theorem}

\begin{Remark}
Since the constant $c_1$ is independent of $\rho_0$, the condition $\|\rho_0-\overline{\rho}\|_{L^2}\le c_1$
allows the initial density  to have large oscillation on the set with small measure.
Given $\rho_0$, if we take the the highly oscillating initial velocity, for example,
\beno
u_0(x)=\sin(\f x \varepsilon)\varphi(x)
\eeno
for any $\varphi\in \cS(\R^3)$, then for $p>3$ and $\varepsilon$ small enough,
\beno
\|u_0\|_{\dot H^{-\delta}}+\|u_0\|_{\dot B^{\f3p-1}_{p,1}}\le C\varepsilon^{\min(\delta,1-\f 3p)}\le c_2.
\eeno
Interestingly, $u_0$ has a large norm in $H^s$ for $s\ge 0$, hence the large energy.
\end{Remark}

\begin{Remark}
Whether similar result can be generalized to the case of two dimensional space remains unknown.
However, it seems impossible to generalize this result to the full Navier-Stokes equations, see \cite{CMZ-Arxiv}.
\end{Remark}

The main steps of the proof are as follows. First of all, we construct a local solution of
(\ref{eq:cNS}) for the initial data $(\rho_0-\overline{\rho}, u_0)$
belonging to the critical Besov space $\dot B^{\f 3p}_{p,1}\times \dot B^{\f3p-1}_{p,1}$.
To deal with the case of physical non vacuum (i.e., $\rho_0\ge c_0>0$),
we need to use the weighted Besov space introduced by Chen, Miao and Zhang \cite{CMZ-Rev}.
Using the smoothing effect of the momentum equation, this solution can propagate the smallness of the initial velocity in Sobolev space of
negative regularity index to that of the velocity at some time $t_0>0$ in Sobolev space of positive regularity index.
Then Hoff's energy method  can be applied to prove that the density is bounded below and above after $t>t_0$.
Finally, we prove a continuation criterion of smooth solution, which ensures that the local solution can be extended to a global solution.

\setcounter{equation}{0}
\section{Littlewood-Paley analysis}

\subsection{Littlewood-Paley decomposition}

Let us introduce the Littlewood-Paley decomposition. Choose a
radial function  $\varphi \in {\cS}(\R^3)$ supported in
${\cC}=\{\xi\in\R^3,\, \frac{3}{4}\le|\xi|\le\frac{8}{3}\}$ such
that
\beno \sum_{j\in\Z}\varphi(2^{-j}\xi)=1 \quad \textrm{for
all}\,\,\xi\neq 0. \eeno The frequency localization operator
$\Delta_j$ and $S_j$ are defined by
\begin{align}
\Delta_jf=\varphi(2^{-j}D)f,\quad S_jf=\sum_{k\le
j-1}\Delta_kf\quad\mbox{for}\quad j\in \Z. \nonumber
\end{align}
With our choice of $\varphi$, one can easily verify that
\beq\label{orth}
\begin{aligned}
&\Delta_j\Delta_kf=0\quad \textrm{if}\quad|j-k|\ge 2\quad
\textrm{and}
\quad \\
&\Delta_j(S_{k-1}f\Delta_k f)=0\quad \textrm{if}\quad|j-k|\ge 5.
\end{aligned}
\eeq

Next we recall Bony's decomposition from \cite{Bony}:
\beq\label{Bonydecom}
uv=T_uv+T_vu+R(u,v), \eeq with
$$T_uv=\sum_{j\in\Z}S_{j-1}u\Delta_jv, \quad R(u,v)=\sum_{j\in\Z}\Delta_ju \widetilde{\Delta}_{j}v,
\quad \widetilde{\Delta}_{j}v=\sum_{|j'-j|\le1}\Delta_{j'}v.$$

The following Bernstein lemma will be frequently used(see \cite{Bah}).

\begin{Lemma}\label{Lem:Bernstein}
Let $1\le p\le q\le+\infty$. Assume that $f\in L^p(\R^3)$, then for
any $\gamma\in(\N\cup\{0\})^3$, there exist constants $C_1$, $C_2$
independent of $f$, $j$ such that
\beno
&&{\rm supp}\hat f\subseteq
\{|\xi|\le A_02^{j}\}\Rightarrow \|\partial^\gamma f\|_q\le
C_12^{j{|\gamma|}+3j(\frac{1}{p}-\frac{1}{q})}\|f\|_{p},
\\
&&{\rm supp}\hat f\subseteq \{A_12^{j}\le|\xi|\le
A_22^{j}\}\Rightarrow \|f\|_{p}\le
C_22^{-j|\gamma|}\sup_{|\beta|=|\gamma|}\|\partial^\beta f\|_p.
\eeno
\end{Lemma}

\subsection{The functional spaces}

We denote the space ${\cZ'}(\R^3)$ by the dual space of
${\cZ}(\R^3)=\{f\in {\cS}(\R^3);\,D^\alpha \hat{f}(0)=0;
\forall\alpha\in\N^3 \,\mbox{multi-index}\}$, which  can be
identified by the quotient space of ${\cS'}(\R^3)/{\cP}$ with the
polynomials space ${\cP}$.

\begin{Definition} Let $s\in\R$, $1\le p,
r\le+\infty$. The homogeneous Besov space $\dot{B}^{s}_{p,r}$ is
defined by
$$\dot{B}^{s}_{p,r}=\{f\in {\cZ'}(\R^3):\,\|f\|_{\dot{B}^{s}_{p,r}}<+\infty\},$$
where \beno
\|f\|_{\dot{B}^{s}_{p,r}}\eqdefa \Bigl\|2^{ks}
\|\Delta_kf(t)\|_{p}\Bigr\|_{\ell^r}.\eeno
In particular, $\dot B^{s}_{2,2}=\dot H^{s}$, $\dot H^s$ is the homogeneous Sobolev space.
\end{Definition}

Let $\{e_k(t)\}_{k\in\Z}$ be a sequence defined in $[0,+\infty)$
satisfying the following conditions
\ben\label{assu:e_k}
e_k(t)\in [0,1], \quad e_k(t)\le e_{k'}(t)\quad \textrm{if}\quad k\le k'\quad \textrm{and} \quad
e_k(t)\sim e_{k'}(t)\quad \textrm{if}\quad k\sim k',
\een
where $k\sim k'$ means that there exists a constant $N_0$ such that $|k-k'|\le N_0$.
The weight function $\{\om_k(t)\}_{k\in\Z}$ is defined by
$$
\om_k(t)=\sum_{\ell\ge k}2^{k-\ell}e_{\ell}(t),\quad k\in \Z.
$$
It is easy to verify that for any $k\in \Z$,
\begin{equation} \label{weightprop}
\begin{split}
&\om_k(t)\le 2,\quad e_k(t)\le \om_k(t), \\
&\om_k(t)\le 2^{k-k'}\om_{k'}(t)\quad \textrm{if}\,\, k\ge k',\quad
\om_k(t)\le 3\om_{k'}(t)\quad \textrm{if}\,\, k\le k',
\\& \om_k(t)\sim \om_{k'}(t)\quad \textrm{if}\,\, k\sim k'.
\end{split}
\end{equation}
The following weighted Besov space is introduced in \cite{CMZ-Rev}.

\begin{Definition}\label{Def:weightspacefuction}
Let $s\in\R$, $1\le p,r\le+\infty$, $0<T<+\infty$. The weighted
Besov space $\dot B^s_{p,r}(\om)$ is  defined by
$$\dot B^s_{p,r}(\om)=\{f\in {\cZ'}(
\R^{N}):\,\|f\|_{\dot
B^s_{p,r}(\om)}<+\infty\},$$ where
\beno
\|f\|_{\dot B^s_{p,r}(\om)}\eqdefa\bigl\|2^{ks}\om_k(T)
\|\Delta_kf\|_{p}\bigr\|_{\ell^r}. \eeno
\end{Definition}
Obviously, $\dot B^s_{p,r}\subset \dot B^s_{p,r}(\om)$ and
\beno
\|f\|_{\dot B^s_{p,r}(\om)}\le 2\|f\|_{\dot B^s_{p,r}}.
\eeno

We also need to use Chemin-Lerner type Besov spaces introduced in \cite{Chemin-JDE}.

\begin{Definition}Let $s\in\R$, $1\le
p,q,r\le+\infty$, $0<T\le+\infty$. The functional space
$\widetilde{L}^q_T(\dot{B}^s_{p,r})$ is defined as the set of all
the distributions $f$ satisfying
$$\|f\|_{\widetilde{L}^q_T(\dot{B}^{s}_{p,r})}\eqdefa \Bigl\|2^{ks}
\|\Delta_kf(t)\|_{L^q(0,T;L^p)}\Bigr\|_{\ell^r}<+\infty.$$
\end{Definition}
The weighted functional space $\widetilde{L}^q_T(\dot B^s_{p,1}(\om))$ is defined similarly,
whose norm is given by
\beno
\|f\|_{\widetilde{L}^q_T(\dot B^s_{p,1}(\om))}\eqdefa\Big\|2^{ks}\om_k(T)
\|\Delta_kf(t)\|_{L^q(0,T;L^p)}\Big\|_{\ell^r}.
\eeno

\subsection{Nonlinear estimates in Besov spaces}

Let us recall some nonlinear estimates in weighted Besov spaces from \cite{CMZ-Rev}.

\begin{Lemma}\label{Lem:bonyweighttimespace}
Let $1\le p,q,q_1,q_2\le \infty$ with $\f 1{q_1}+\f1{q_2}=\f1q$.
Then there hold

\no (a)\, if $s_2\le \frac{3}{p}$, we have
\beno
\|T_gf\|_{\widetilde{L}^{q}_T(\dot B^{s_1+s_2-
\frac{3}{p}}_{p,1}(\omega))} \le C\|f\|_{\widetilde{L}^{q_1}_T(\dot
B^{s_1}_{p,1}(\omega))}\|g\|_{\widetilde{L}^{q_2}_T(\dot
B^{s_2}_{p,1})}; \eeno

\no (b)\, if $s_1\le \frac{3}{p}-1$, we have
\beno
\|T_fg\|_{\widetilde{L}^{q}_T(\dot B^{s_1+s_2-
\frac{3}{p}}_{p,1}(\omega))} \le C\|f\|_{\widetilde{L}^{q_1}_T(\dot
B^{s_1}_{p,1}(\omega))}\|g\|_{\widetilde{L}^{q_2}_T(\dot
B^{s_2}_{p,1})};
\eeno

\no (c)\, if $s_1+s_2>3\max(0,\frac 2p-1)$, we have
\beno
\|R(f,g)\|_{\widetilde{L}^{q}_T(\dot B^{s_1+s_2-
\frac{3}{p}}_{p,1}(\omega))} \le C\|f\|_{\widetilde{L}^{q_1}_T(\dot
B^{s_1}_{p,1}(\omega))}\|g\|_{\widetilde{L}^{q_2}_T(\dot
B^{s_2}_{p,1})}.
\eeno
\end{Lemma}

The following lemma is a direct consequence of Lemma \ref{Lem:bonyweighttimespace}.

\begin{Lemma}\label{Lem:bi-weight}
Let $s_1\le \frac{3}{p}-1, s_2\le \frac{3}{p}, s_1+s_2>3\max(0,\frac 2p-1)$, and $1\le p,q,q_1,q_2\le \infty$
with $\f 1{q_1}+\f1{q_2}=\f1q$. Then there holds
\beno
\|fg\|_{\widetilde{L}^{q}_T(\dot B^{s_1+s_2-
\frac{3}{p}}_{p,1}(\omega))} \le C\|f\|_{\widetilde{L}^{q_1}_T(\dot
B^{s_1}_{p,1}(\omega))}\|g\|_{\widetilde{L}^{q_2}_T(\dot
B^{s_2}_{p,1})}.
\eeno
\end{Lemma}

In general Besov spaces, we have

\begin{Lemma}\label{Lem:bi-est}
Let $s_1, s_2\le \frac{3}{p},\, s_1+s_2>3\max (0,\frac2p-1)$, and $1\le p,q,q_1,q_2\le \infty$ with $\f 1{q_1}+\f1{q_2}=\f1q$.
Then there holds
\beno \|fg\|_{\widetilde{L}^{q}_T(\dot B^{s_1+s_2-
\frac{3}{p}}_{p,1})} \le C\|f\|_{\widetilde{L}^{q_1}_T(\dot
B^{s_1}_{p,1})}\|g\|_{\widetilde{L}^{q_2}_T(\dot B^{s_2}_{p,1})}.
\eeno
\end{Lemma}

\begin{Lemma}\label{Lem:non-est}
Let $s>0$ and $1\le p,q,r\le \infty$. Assume that $F\in W^{[s]+3,\infty}_{loc}(\R)$ with  $F(0)=0$.
Then there holds
\beno
\|F(f)\|_{\widetilde{L}^q_T(\dot B^s_{p,r}(\om))}
\le C(1+\|f\|_{L^\infty_T(L^\infty)})^{[s]+2}\|f\|_{\widetilde{L}^q_T(\dot B^s_{p,r}(\om))}.
\eeno
The same result also holds true in Besov space without weight.
\end{Lemma}

The following is a weighted commutator estimate.

\begin{Lemma}\label{Lem:com-est}
Let $p\in [1,\infty)$ and $s\in (-3\min(\frac 1p,\frac1{p'}),\frac 3p]$.
Then  there holds
\beno
\big\|2^{js}\omega_j(T)\|[\Delta_j,f]\na
g\|_{L^1_T(L^p)}\big\|_{\ell^1}\le C\|f\|_{\widetilde{L}^\infty_T(\dot B^{
\frac{3}{p}}_{p,1}(\om))}\|g\|_{\widetilde{L}^1_T(\dot
B^{s+1}_{p,1})}.
\eeno
\end{Lemma}

\no{\bf Proof.}\,Using Bony's decomposition (\ref{Bonydecom}), we write
\beno
[f,\Delta_j]\cdot\p_kg&=&[T_{f},\Delta_j]\p_kg+T_{\p_k\Delta_jg}f+R(f,\p_k\Delta_jg)\\
&&-\Delta_j(T_{\p_kg}f)-\Delta_jR(f,\p_kg).
\eeno
Using Lemma \ref{Lem:bonyweighttimespace} (a) and (c) with $s_1= \frac{3}{p}$ and
$s_2=s$, we get
\beno
&&\big\|\om_j(T)2^{js}\|\Delta_j(T_{\p_kg}f)\|_{L^1_T(L^p)}\big\|_{\ell^1}
\le C\|f\|_{\widetilde{L}^\infty_T(\dot B^{ \frac{3}{p}}_{p,1}(\om))}\|g\|_{\widetilde{L}^1_T(\dot B^{s+1}_{p,1})},\\
&&\big\|\om_j(T)2^{js}\|\Delta_jR(f^k,\p_kg)\|_{L^1_T(L^p)}\big\|_{\ell^1}\le
C\|f\|_{\widetilde{L}^\infty_T(\dot B^{
\frac{3}{p}}_{p,1}(\om))}\|g\|_{\widetilde{L}^1_T(\dot
B^{s+1}_{p,1})}.
\eeno
Noticing that
\ben\label{eq:para}
T_{\p_k\Delta_jg}'f\triangleq T_{\p_k\Delta_jg}f+R(f,\p_k\Delta_jg)=\sum_{j'\ge
j-2}S_{j'+2}\Delta_j\p_kg\Delta_{j'}f,
\een
then we get by Lemma \ref{Lem:Bernstein} and (\ref{weightprop}) that
\beno
&&\big\|\om_j(T)2^{js}\|T_{\p_k\Delta_jg}'f\|_{L^1_T(L^p)}\big\|_{\ell^1}\\
&&\le C\big\|\om_j(T)2^{js}\|\Delta_j\na g\|_{L^1_T(L^\infty)}\sum_{j'\ge j-2}\|\Delta_{j'}f\|_{L^\infty_T(L^p)}\big\|_{\ell^1}\\
&&\le C\big\|2^{j(s+1+\frac{3}{p})}\|\Delta_jg\|_{L^1_T(L^p)}\sum_{j'\ge j-2}\om_{j'}(T)\|\Delta_{j'}f\|_{L^\infty_T(L^p)}\big\|_{\ell^1}\\
&&\le C\|f\|_{\widetilde{L}^\infty_T(\dot B^{
\frac{3}{p}}_{p,1}(\om))}\|g\|_{\widetilde{L}^1_T(\dot
B^{s+1}_{p,1})}.
\eeno
Set $h(x)=(\cF^{-1}\varphi)(x)$. Then we have
\beno [T_{f},
\Delta_j]\pa_k g
&=&\sum_{|j'-j|\le4}2^{4j}\int_{\R^3}\int_0^1y\cdot\na
S_{j'-1}f(x-\tau
y)d\tau\pa_kh(2^jy)\Delta_{j'}g(x-y)dy\nonumber\\&&\qquad\quad+
2^{3j}\int_{\R^3}h(2^j(x-y))\pa_kS_{j'-1}f(y)\Delta_{j'}g(y)dy ,
\eeno
from which and Young's inequality, we infer that
\beno
\big\|\om_j(T)2^{js}\|[T_{f}, \Delta_j]\pa_k
g\|_{L^1_T(L^p)}\big\|_{\ell^1}\le C\|f\|_{\widetilde{L}^\infty_T(\dot B^{
\frac{3}{p}}_{p,1}(\om))}\|g\|_{\widetilde{L}^1_T(\dot
B^{s+1}_{p,1})}.
\eeno

The proof is finished by summing up all the above estimates.\ef

\setcounter{equation}{0}
\section{The linear transport equation and momentum equation}

We consider the linear transport equation
\begin{align}\label{eq:linertrans}
\bigg\{\begin{aligned}
&\partial_t f+v\cdot \nabla f =g,\\
&f(0,x)=f_0.
\end{aligned}
\bigg.\end{align}

Set $V(t)\eqdefa\int_0^t\|\nabla
v(\tau)\|_{\dot{B}^{\frac{3}{p}}_{p,1}}d\tau.$ The following result
is from \cite{Dan-Non}.

\begin{Proposition}\label{Prop:transport}
Let $p\in [1,+\infty]$ and $s\in (-3\min(\frac1p,\frac1{p'}), 1+\frac 3p)$.
Let $v$ be a vector field such that $\nabla v\in L^1_T(\dot{B}^{\frac{3}{p}}_{p,1})$. Assume that
$f_0\in \dot{B}^{s}_{p,1},$ $g\in L^1_T(\dot{B}^{s}_{p,1})$ and $f$
is the solution of (\ref{eq:linertrans}). Then there holds for
$t\in[0,T]$, \beno
\|f\|_{\widetilde{L}^\infty_t(\dot{B}^{s}_{p,1})}\le e^{CV(t)}\Big(
\|f_0\|_{\dot{B}^{s}_{p,1}}+\int_0^t
e^{-CV(\tau)}\|g(\tau)\|_{\dot{B}^{s}_{p,1}}d\tau\Big).
\eeno

\end{Proposition}

The following weighted version is from \cite{CMZ-Rev}.

\begin{Proposition}\label{Prop:transportweight}
Let $p\in [1,+\infty]$ and $s\in (-3\min(\frac1p,\frac1{p'}),\frac
3p]$. Let $v$ be a vector field such that $\nabla v\in
L^1_T(\dot{B}^{ \frac{3}{p}}_{p,1})$. Assume that $f_0\in
\dot{B}^{s}_{p,1},$ $g\in L^1_T(\dot{B}^{s}_{p,1})$ and $f$ is the
solution of (\ref{eq:linertrans}). Then there holds for $t\in
[0,T]$,
\beno
\|f\|_{\widetilde{L}^\infty_t(\dot
B^s_{p,1}(\omega))}\le e^{CV(t)}\Big(
\|f_0\|_{\dot{B}^{s}_{p,1}(\om)}+\int_0^t
e^{-CV(\tau)}\|g(\tau)\|_{\dot{B}^{s}_{p,1}(\om)}d\tau\Big).
\eeno
\end{Proposition}

We next study the linear momentum equations with variable coefficients
\ben\label{eq:linearmomen} \left\{
\begin{array}{ll}
\p_tu-\dv(\overline{\mu}\na u)
-\na(\overline{\lambda}\dv\,u)=G, \\
u|_{t=0}=u_0.
\end{array}
\right.
\een
Assume that the viscosity coefficients $\overline{\lambda}(\rho)$ and $\overline{\mu}(\rho)$ depend smoothly on the function $\rho$
and there exists a positive constant $c_3$ such that
\beno
\overline{\mu}\ge c_3,\quad \overline{\lambda}+\overline{\mu}\ge c_3.
\eeno
And the weighted function $\om_k(t)$ is taken as follows
$$
\om_k(t)=\sum_{\ell\ge k}2^{k-\ell}e_\ell(t),
$$
with $e_\ell(t)=(1-e^{-c2^{2\ell}t})^\f12$ for some positive constant $c>0$.
It is easy to verify that the function $e_\ell(t)$ satisfies (\ref{assu:e_k}).

\begin{Proposition}\label{Prop:momenequ}
Let $p\in [2,\infty), s\in (1-\f 3p,\frac 3p], q\in [1,\infty]$.
Assume that $u_0\in \dot B^{s-1}_{p,1}, G\in L^1_T(\dot B^{s-1}_{p,1})$,
and $\rho-\overline{\rho}\in \widetilde{L}^\infty_T(\dot B^{
\frac{3}{p}}_{p,1})$. Let $u$ be a solution of (\ref{eq:linearmomen}).
Then we have
\beno
&&\|u\|_{\widetilde{L}^1_T(\dot
B^{s+1}_{p,1})}+\|u\|_{\widetilde{L}^2_T(\dot B^{s}_{p,1})}\\
&& \le C\Bigl(\|u_0\|_{\dot B^{s-1}_{p,1}(\om))}
+\|G\|_{\widetilde{L}^1_T(\dot B^{s-1}_{p,1}(\om))}
+A(T)\|\rho-\overline{\rho}\|_{\widetilde{L}^\infty_T(\dot B^{
\frac{3}{p}}_{p,1}(\om))}\|u\|_{\widetilde{L}^1_T(\dot
B^{s+1}_{p,1})}\Bigr),
\eeno
and the version without weight
\beno
\|u\|_{\widetilde{L}^q_T(\dot B^{s-1+\f 2q}_{p,1})}\le C\Bigl(\|u_0\|_{\dot
B^{s-1}_{p,1}}+\|G\|_{\widetilde{L}^1_T(\dot B^{s-1}_{p,1})}
+A(T)\|\rho-\overline{\rho}\|_{\widetilde{L}^\infty_T(\dot B^{
\frac{3}{p}}_{p,1})}\|u\|_{\widetilde{L}^1_T(\dot
B^{s+1}_{p,1})}\Bigr).
\eeno
Here $A(T)\eqdefa\bigl(1+\|\rho\|_{L^\infty_T(L^\infty)}\bigr)^{[\frac{3}{p}]+2}$.
\end{Proposition}

\no{\bf Proof.}\, We only prove the first inequality,
the proof of the second inequality is similar.
Set $d=\dv u$ and $w=\textrm{curl} u$.
By (\ref{eq:linearmomen}), $(d,w)$ satisfies
\ben\label{eq:dv-curl} \left\{
\begin{array}{ll}
\p_td-\dv(\overline{\nu}\na d)=\dv G+F_1, \\
\p_tw-\dv(\overline{\mu}\na w)=\textrm{curl} G+F_2,\\
(d,w)|_{t=0}=(\dv u_0,\textrm{ curl} u_0)\triangleq (d_0,w_0),
\end{array}
\right. \een
where $\overline{\nu}=\overline{\lambda}+\overline{\mu}$ and
\beno
&&F_1=\dv(\na \overline{\mu}\cdot\na u)+\dv(\na (\overline{\lambda}+\overline{\mu})d),\\
&&F_2^{i,j}=\dv(\p_j\overline{\mu}\na u^i-\p_i\overline{\mu}\na u^j),\quad i,j=1,2,3.
\eeno
Apply the operator $\Delta_j$ to (\ref{eq:dv-curl}) to obtain
\beno \left\{
\begin{array}{ll}
\p_t\Delta_jd-\dv(\overline{\nu}\na \Delta_jd)=\dv \Delta_jG+\Delta_jF_1+\dv[\Delta_j,\overline{\nu}]\na d, \\
\p_t\Delta_jw-\dv(\overline{\mu}\na \Delta_jw)=\textrm{curl} \Delta_jG+\Delta_jF_2
+\dv[\Delta_j,\overline{\mu}]\na w.
\end{array}
\right. \eeno
Multiplying the first equation by $|\Delta_jd|^{p-2}\Delta_jd$, we get by integrating over $\R^3$ that
\beno
&&\f1p\f d{dt}\|\Delta_jd\|_p^p-\int_{\R^3}\dv(\overline{\nu}\na \Delta_jd)|\Delta_jd|^{p-2}\Delta_jd \intd x\\
&&\quad=\int_{\R^3}\bigl(\dv
\Delta_jG+\Delta_jF_1+\dv[\Delta_j,\overline{\nu}]\na
d\bigr)|\Delta_jd|^{p-2}\Delta_jd dx.
\eeno
There exists $c_p>0$ depending on $c_3, p$ such that(see \cite{Dan-CPDE})
\beno
\f1p\f d{dt}\|\Delta_jd\|_p^p+c_p2^{2j}\|\Delta_jd\|_p^p\le
C\|\Delta_jd\|_p^{p-1} \bigl(2^j\|
\Delta_jG\|_p+\|\Delta_jF_1\|_p+2^j\|[\Delta_j,\overline{\nu}]\na
d\|_p\bigr). \eeno
This gives that
\beno
\f d{dt}\|\Delta_jd\|_p+c_p2^{2j}\|\Delta_jd\|_p\le C\big(2^j\|\Delta_jG\|_p
+\|\Delta_jF_1\|_p+2^j\|[\Delta_j,\overline{\nu}]\na d\|_p\big),
\eeno
which implies that
\beno
\|\Delta_jd(t)\|_p\le e^{-c_p2^{2j}t}\|\Delta_jd_0\|_p+C\int_0^te^{-c_p2^{2j}(t-\tau)}\bigl(2^j\|
\Delta_jG\|_p+\|\Delta_jF_1\|_p+2^j\|[\Delta_j,\overline{\nu}]\na
d\|_p\bigr)d\tau.
\eeno
Similarly, we can show that
\beno
\|\Delta_jw(t)\|_p\le
e^{-c_p2^{2j}t}\|\Delta_jd_0\|_p
+C\int_0^te^{-c_p2^{2j}(t-\tau)}\bigl(2^j\|\Delta_jG\|_p+\|\Delta_jF_2\|_p+2^j\|[\Delta_j,\overline{\mu}]\na
w\|_p\bigr)d\tau.
\eeno
Hence, we infer that for any $q\in [1,\infty]$ and $t\in [0,T]$,
\ben\label{eq:localcurl-dvest}
&&\|\Delta_jd(t)\|_{L^q_t(L^p)}+\|\Delta_jw(t)\|_{L^q_t(L^p)}\nonumber\\
&&\le C2^{-2j/q}c_j(T)^\f1q(\|\Delta_jd_0\|_p+\|\Delta_jw_0\|_p)\nonumber\\
&&\quad+C2^{-2j/q}c_j(T)^\f1q\bigl(2^j\|\Delta_jG\|_{L^1_t(L^p)}+\|\Delta_jF_1\|_{L^1_t(L^p)}+
2^j\|[\Delta_j,\overline{\nu}]\na d\|_{L^1_t(L^p)}\bigr)\nonumber\\
&&\quad+C2^{-2j/q}c_j(T)^\f1q\bigl(2^j\|\Delta_jG\|_{L^1_t(L^p)}+\|\Delta_jF_2\|_{L^1_t(L^p)}
+2^j\|[\Delta_j,\overline{\mu}]\na w\|_{L^1_t(L^p)}\bigr),
\een
with $c_j(T)=1-e^{-c_p2^{2j}T}$. Notice that
\beno
2^j\|\Delta_ju\|_p\sim \|\Delta_jd\|_p+\|\Delta_jw\|_p, \quad e_j(T)\le \om_j(T),
\eeno
which along with (\ref{eq:localcurl-dvest}) implies that
\ben\label{eq:momenone}
&&\|u\|_{\widetilde{L}^1_T(\dot B^{s+1}_{p,1})}+\|u\|_{\widetilde{L}^2_T(\dot B^{s}_{p,1})}\nonumber\\
&&\le C\bigl(\|u_0\|_{\dot B^{s-1}_{p,1}(\om)}+\|G\|_{\widetilde{L}^1_T(\dot B^{s-1}_{p,1}(\om))}+\|(F_1,F_2)\|_{\widetilde{L}^1_T(\dot B^{s-2}_{p,1}(\om))}\bigr)\nonumber\\
&&\qquad+C\big\|2^{j(s-1)}w_j(T)\bigl(\|[\Delta_j,\overline{\nu}]\na
d\|_{L^1_T(L^p)}+\|[\Delta_j,\overline{\mu}]\na
w\|_{L^1_T(L^p)}\bigr)\big\|_{\ell^1}.
\een
From Lemma \ref{Lem:bi-weight} and Lemma \ref{Lem:non-est},
it follows that
\ben\label{eq:F1-w}
\|F_1\|_{\widetilde{L}^1_T(\dot B^{s-2}_{p,1}(\om))}+\|F_2\|_{\widetilde{L}^1_T(\dot B^{s-2}_{p,1}(\om))}
\le A(T)\|\rho-\overline{\rho}\|_{\widetilde{L}^\infty_T(\dot
B^{ \frac{3}{p}}_{p,1}(\om))}\|u\|_{\widetilde{L}^1_T(\dot
B^{s+1}_{p,1})},
\een
and using Lemma \ref{Lem:com-est}, we get
\ben\label{eq:comm-w}
&&\big\|2^{j(s-1)}\omega_j(T)\bigl(\|[\Delta_j,\overline{\nu}]\na d\|_{L^1_T(L^p)}+\|[\Delta_j,\overline{\mu}]\na w\|_{L^1_T(L^p)}\bigr)\big\|_{\ell^1}\nonumber\\
&&\le A(T)\|\rho-\overline{\rho}\|_{\widetilde{L}^\infty_T(\dot B^{
\frac{3}{p}}_{p,1}(\om))}\|u\|_{\widetilde{L}^1_T(\dot
B^{s+1}_{p,1})}.
\een
Summing up (\ref{eq:momenone})-(\ref{eq:comm-w}) yields the first inequality
of Proposition \ref{Prop:momenequ}. \ef

\setcounter{equation}{0}

\section{Local well-posedness with physical non vacuum}

In this section, we prove the local well-posedness of compressible Navier-Stokes equations
in the critical Besov space with physical non vacuum assumption.
That is, we just impose $\rho_0\ge c_00$ for the initial density.
In the case when $\rho-\overline{\rho}_0$ is small in $\dot B^{ \frac{3}{p}}_{p,1}$ or has more regularity,
the corresponding local-well posedness has been obtained by Danchin \cite{Dan-CPDE}.
Recently, Chen-Miao-Zhang \cite{CMZ-Rev} and Danchin \cite{Dan-CPDE07} developed two different methods to deal with physical non vacuum case.
Here we will revisit the proof of \cite{CMZ-Rev}. The main goal is to record the precise dependence of energy
bounds and the existence time $T$ on the initial data and other important constants, which is very important for our argument.

\begin{Theorem}\label{thm:local}
Let $\overline{\rho}>0$ and $c_0>0$. Assume that the initial data
$(\rho_0,u_0)$ satisfies
\ben
\rho_0-\overline{\rho}\in \dot B^{\f 3 p}_{p,1},\quad c_0\le \rho_0\le c_0^{-1}, \quad u_0\in \dot B^{\f 3p-1}_{p,1}.
\een
Then there exists a positive time $T$ such that\vskip 0.1cm

(a)\, {\bf Existence:}\, If $p\in [2,6)$, the system (\ref{eq:cNS}) has a solution $(\rho-\bar{\rho},u)\in E^p_T$ with
\beno
E^p_T\eqdefa C([0,T]; \dot B^{\f {3} p}_{p,1})\times \Bigl(C([0,T]; \dot B^{\f {3}
p-1}_{p,1}) \cap L^1(0,T; \dot B^{\f {3} p+1}_{p,1})\Bigr)^3,\quad \rho\ge
\f12 c_0;
\eeno

(b)\, {\bf Uniqueness:} If $p\in (1,3]$, then the uniqueness holds in $E^p_T$.\vskip 0.1cm
\end{Theorem}

\no{\bf Proof.}\,We only present the uniform energy estimates for smooth solution
of (\ref{eq:cNS}) , since the existence can be deduced by constructing approximate solution sequence and
compact argument. The uniqueness has been proved in \cite{CMZ-Rev}.
Set
\beno a(t,x)=\f {\rho(t,x)-\overline{\rho}}
{\overline{\rho}},\quad \overline{\mu}(\rho)=\f
{\mu}{\rho},\quad
\overline{\lambda}(\rho)=\f
{\lambda}{\rho}.
\eeno
Then the system (\ref{eq:cNS}) can be rewritten as
\ben\label{eq:linearcNS} \left\{
\begin{array}{ll}
\p_ta+u\cdot\na a=F,\\
\p_tu-\dv(\overline{\mu}\na u)
-\na(\overline{\lambda}\dv\,u)=G, \\
(a,u)|_{t=0}=(a_0,u_0),
\end{array}
\right.
\een
with $a_0=\f {\rho_0(x)-\bar{\rho}} {\bar{\rho}}$ and
\beno
&&F(=-(1+a)\dv\, u,\\
&&G=-u\cdot\na u-\f{\overline{\rho}P'(\rho)} {\rho}\na a+\f
{\mu} {\rho^2}\na \rho\cdot\na u+\f
{\lambda}{\rho^2}\na \rho\dv \,u.
\eeno

Let us assume that the solution $(\rho,u)$ satisfies
\ben
&&\f {c_0} 2\le \rho(t,x)\le 2c_0^{-1}\quad\textrm{ on }\quad [0,T];\label{eq:H1-a}\\
&&\|a\|_{\widetilde{L}^\infty_T(\dot B^{
\frac{3}{p}}_{p,1})}+\|u\|_{\widetilde{L}^\infty_T(\dot B^{
\frac{3}{p}-1}_{p,1})}\le C_0\big(\|a_0\|_{\dot B^{ \frac{3}{p}}_{p,1}}+\|u_0\|_{\dot B^{
\frac{3}{p}-1}_{p,1}}\big)\eqdefa C_0E_0;\label{eq:H2-av}\\
&&\|a\|_{\widetilde{L}^\infty_T(\dot B^{
\frac{3}{p}}_{p,1}(\om))}\le \eta_1;\label{eq:H3-a}\\
&&\|u\|_{\widetilde{L}^1_T(\dot B^{ \frac{3}{p}}_{p,1})}
+\|u\|_{\widetilde{L}^2_T(\dot B^{ \frac{3}{p}}_{p,1})}\le
\eta_2.\label{eq:H4-v}
\een
These assumptions are satisfied for the solution in $E^p_T$ if $T$ is small enough.
We will show that there exists a $T>0$ depending on the initial data and $\lambda,\mu,\overline{\rho}, c_0, p$ such that these assumptions are satisfied.

In what follows, we denote $C$ by a constant depending only on $\lambda,\mu,\overline{\rho}, c_0, p$,
which may be different from line to line.

{\bf Step 1.\,} Estimates in the Besov space

 Applying Proposition \ref{Prop:transport} to the first equation of (\ref{eq:linearcNS}), we get
\ben\label{eq:dens1}
\|a\|_{\widetilde{L}^\infty_T(\dot{B}^{\frac{3}{p}}_{p,1})}
\le e^{CV(T)}\bigl( \|a_0\|_{\dot{B}^{
\frac{3}{p}}_{p,1}}+\|F\|_{\widetilde{L}^1_T(\dot{B}^{
\frac{3}{p}}_{p,1})}\bigr),
\een
and applying Proposition \ref{Prop:momenequ} to the second equation of (\ref{eq:linearcNS}),
we have
\ben\label{eq:veloc1}
\|u\|_{\widetilde{L}^\infty_T(\dot B^{\frac{3}{p}-1}_{p,1})}
\le C\|u_0\|_{\dot B^{ \frac{3}{p}-1}_{p,1}}+C\|G\|_{\widetilde{L}^1_T(\dot B^{ \frac{3}{p}-1}_{p,1})}+C\|a\|_{\widetilde{L}^\infty_T(\dot B^{
\frac{3}{p}}_{p,1})}\|u\|_{\widetilde{L}^1_T(\dot B^{
\frac{3}{p}+1}_{p,1})},
\een
where $V(t)=\int_0^t\|u(\tau)\|_{\dot{B}^{\frac{3}{p}+1}_{p,1}}d\tau$.
By Lemma \ref{Lem:bi-est}, we have
\beno
\|F\|_{\widetilde{L}^1_T(\dot{B}^{ \frac{3}{p}}_{p,1})}\le
\|u\|_{\widetilde{L}^1_T(\dot{B}^{ \frac{3}{p}+1}_{p,1})}+
C\|a\|_{\widetilde{L}^\infty_T\dot{B}^{
\frac{3}{p}}_{p,1}}\|u\|_{\widetilde{L}^1_T(\dot{B}^{
\frac{3}{p}+1}_{p,1})},
\eeno
and by Lemma \ref{Lem:bi-est} and Lemma \ref{Lem:non-est},
\beno
\|G\|_{\widetilde{L}^1_T(\dot B^{ \frac{3}{p}-1}_{p,1})}
\le C\|u\|_{\widetilde{L}^\infty_T(\dot{B}^{ \frac{3}{p}-1}_{p,1})}\|u\|_{\widetilde{L}^1_T(\dot{B}^{ \frac{3}{p}+1}_{p,1})}
+C\|a\|_{\widetilde{L}^\infty_T(\dot{B}^{
\frac{3}{p}}_{p,1})} \bigl(T+\|u\|_{\widetilde{L}^1_T(\dot{B}^{
\frac{3}{p}+1}_{p,1})}\bigr).
\eeno
Here we used $p<6$.
Plugging them into (\ref{eq:dens1})-(\ref{eq:veloc1}), we obtain
\ben\label{eq:denvel-infty}
\|a\|_{\widetilde{L}^\infty_T(\dot{B}^{\frac{3}{p}}_{p,1})}
+\|u\|_{\widetilde{L}^\infty_T(\dot B^{\frac{3}{p}-1}_{p,1})}
\le Ce^{C\eta_2}\big(E_0+C_0E_0(\eta_2+T)+\eta_2\big).
\een

{\bf Step 2.\,} Estimates in the weighted Besov space

Applying Proposition \ref{Prop:transportweight} to the first equation of (\ref{eq:linearcNS}),
we get
\ben\label{eq:dens2}
\|a\|_{\widetilde{L}^\infty_T(\dot{B}^{
\frac{3}{p}}_{p,1}(\om))}\le e^{CV(T)}\bigl(\|a_0\|_{\dot{B}^{
\frac{3}{p}}_{p,1}(\om)}+\|F\|_{\widetilde{L}^1_T(\dot{B}^{
\frac{3}{p}}_{p,1}(\om))}\bigr),
\een
and by Proposition \ref{Prop:momenequ}, we have
\ben\label{eq:veloc2}
&&\|u\|_{\widetilde{L}^1_T(\dot
B^{\frac{3}{p}+1}_{p,1})}+\|u\|_{\widetilde{L}^2_T(\dot
B^{\frac{3}{p}}_{p,1})}\nonumber\\&&
\le C\Bigl(\|u_0\|_{\dot B^{\frac{3}{p}-1}_{p,1}(\om))} +\|G\|_{\widetilde{L}^1_T(\dot
B^{\frac{3}{p}-1}_{p,1}(\om))}
+\|a\|_{\widetilde{L}^\infty_T(\dot B^{
\frac{3}{p}}_{p,1}(\om))}\|u\|_{\widetilde{L}^1_T(\dot
B^{\frac{3}{p}+1}_{p,1})}\Bigr).
\een
By Lemma \ref{Lem:bi-est}, we have
\beno
\|F\|_{\widetilde{L}^1_T(\dot{B}^{ \frac{3}{p}}_{p,1}(\om))}\le
2\|F\|_{\widetilde{L}^1_T(\dot{B}^{ \frac{3}{p}}_{p,1})}\le C(1+
\|a\|_{\widetilde{L}^\infty_T(\dot{B}^{
\frac{3}{p}}_{p,1})})\|u\|_{\widetilde{L}^1_T(\dot{B}^{
\frac{3}{p}+1}_{p,1})},
\eeno
and by Lemma \ref{Lem:bi-weight} and Lemma \ref{Lem:non-est},
\beno
\|G\|_{\widetilde{L}^1_T(\dot B^{ \frac{3}{p}-1}_{p,1}(\om))} \le
C\|u\|_{\widetilde{L}^2_T(\dot{B}^{ \frac{3}{p}}_{p,1})}^2
+C\|a\|_{\widetilde{L}^\infty_T(\dot{B}^{\frac{3}{p}}_{p,1}(\om))}
\bigl(T+\|u\|_{\widetilde{L}^1_T(\dot{B}^{
\frac{3}{p}+1}_{p,1})}\bigr).
\eeno
Plugging them into (\ref{eq:dens2})-(\ref{eq:veloc2}) yields that
\ben\label{eq:dens-weight}
&&\|a\|_{\widetilde{L}^\infty_T(\dot{B}^{
\frac{3}{p}}_{p,1}(\om))}\le Ce^{C_1\eta_2}\bigl(
\|a_0\|_{\dot{B}^{ \frac{3}{p}}_{p,1}(\om)}+(1+C_0E_0)\eta_2\bigr),\\
&&\|u\|_{\widetilde{L}^1_T(\dot B^{\frac
3p+1}_{p,1})}+\|u\|_{\widetilde{L}^2_T(\dot B^{\frac 3p}_{p,1})}
\le C\big(\|u_0\|_{\dot B^{\frac 3p-1}_{p,1}(\om))}
+\eta^2_2+\eta_1(T+\eta_2)
\big).\label{eq:veloc-onetwo}
\een

{\bf Step 3.\,} Estimate of the density

Let $X(t,x)$ be a solution of
\beno
\f d {dt}
X(t,x)=u(t,X(t,x)),\quad X(0,x)=x,
\eeno
and
we denote by $X^{-1}(t,x)$ the inverse of $X(t,x)$.
Then $a(t,x)$ can be solved as
\beno
a(t,x)=a_0(X^{-1}(t,x))+\int_0^tF(\tau,X(\tau,X^{-1}(t,x)))d\tau,
\eeno
this gives
\ben\label{eq:rho}
\rho(t,x)=\rho_0(X^{-1}(t,x))+\bar{\rho}_0\int_0^tF(\tau,X(\tau,X^{-1}(t,x)))d\tau.
\een
Since we have
\beno
\|F\|_{L^1_T(L^\infty)}\le
\|\na u\|_{L^1_T(L^\infty)}(1+\|a\|_{L^\infty_T(L^\infty)})\le C\eta_2,
\eeno
from which and (\ref{eq:rho}), it follows that
\ben\label{eq:density-est}
\rho(t,x)\ge c_0-C\eta_2,\quad \rho(t,x)\le c_0^{-1}+C\eta_2.
\een

{\bf Step 4.} Continuity argument

Let $C_0=10C$. We first take $\eta_2$ and $\eta_1$ small enough such that
\ben\label{assu:eta}
\begin{split}
&C\eta_1\le \f 18,\quad  e^{C\eta_2}\le 2,\quad C_0\eta_2\le 1,\quad \eta_2\le E_0,\\
&C\eta_2\le \f {c_0}4,\quad 2C(1+C_0E_0)\eta_2 \le \f {\eta_1} 4.
\end{split}
\een
Then we take $T$ small enough such that
\ben\label{assu:time}
C_0T\le 1,\quad T\le \eta_2,\quad\|a_0\|_{\dot{B}^{\frac{3}{p}}_{p,1}(\om)}\le \f {\eta_1} {4C},
\quad \|u_0\|_{\dot{B}^{ \frac{3}{p}-1}_{p,1}(\om)}\le \f {\eta_2} {4C}.
\een
With this choice of $\eta_1,\eta_2$ and $T$, it follows from (\ref{eq:denvel-infty}), (\ref{eq:dens-weight}), (\ref{eq:veloc-onetwo}) and (\ref{eq:density-est}) that
\beno
&&\|a\|_{\widetilde{L}^\infty_T(\dot{B}^{
\frac{3}{p}}_{p,1})}+\|u\|_{\widetilde{L}^\infty_T(\dot B^{\frac
3p-1}_{p,1})} \le \f 45C_0E_0,\\
&&\|a\|_{\widetilde{L}^\infty_T(\dot B^{
\frac{3}{p}}_{p,1}(\om))}\le \f34\eta_1,\\
&&\|u\|_{\widetilde{L}^1_T(\dot B^{ \frac{3}{p}+1}_{p,1})}
+\|u\|_{\widetilde{L}^2_T(\dot B^{ \frac{3}{p}}_{p,1})}\le \f3 4\eta_2,\\
&&\f 34 c_0 \le \rho(t,x)\le \f 54c_0^{-1},
\eeno
which implies that the assumptions (\ref{eq:H1-a})-(\ref{eq:H4-v}) are satisfied
for $T$ defined by (\ref{assu:time}).\ef

\begin{Remark}
In general, the choice of time $T$ depends on $C_0=C_0(\lambda,\mu, \overline{\rho}, c_0,p)$ and the behaviour of the initial data $(\rho_0,u_0)$.
If $\|u_0\|_{\dot{B}^{ \frac{3}{p}-1}_{p,1}}\le c$ and  $a_0\in \dot B^{\f 3p}_{p,1}\cap H^s$ for $s>\f 32$
and $c=c(\lambda,\mu, \overline{\rho}, c_0,p)$ small enough,
then $T$ can be taken as
\ben\label{eq:T-form}
T=\f 1 {C\big(1+\|a_0\|_{\dot B^{\f 3p}_{p,1}\cap H^s}\big)^{(s+\f12)/(s-\f32)}},\quad C=C(\lambda,\mu, \overline{\rho}, c_0,p,s).
\een
Indeed, due to (\ref{assu:eta}), we can take $\eta_1,\eta_2$ as follows
\beno
\eta_1=\f 1 {\widetilde{C}},\quad \eta_2=\f 1 {\widetilde{C}\big(1+\|a_0\|_{\dot B^{\f 3p}_{p,1}}\big)}.
\eeno
for some constant $\widetilde{C}$ depending only on $\lambda,\mu, \overline{\rho}, c_0,p$.
And by Lemma \ref{Lem:Bernstein}, for any $N\in\N$,
\beno
\sum_{j>N}2^{j\f 3p}\|\Delta_j a_0\|_{p}\le C\sum_{j>N}2^{j\f 32}\|\Delta_j a_0\|_{2}\le C2^{-(s-\f32)N}\|a_0\|_{H^s}.
\eeno
And by (\ref{weightprop}), we have
\beno
\sum_{j\le N}2^{j\f 3p}\om_j(T)\|\Delta_j a_0\|_{p}\le \om_N(T)\|a_0\|_{\dot B^{\f 3p}_{p,1}}.
\eeno
Then (\ref{assu:time}) is ensured by taking $T$ as (\ref{eq:T-form}).

\end{Remark}

\setcounter{equation}{0}

\section{Propagation of Sobolev regularity}

\subsection{Propagation of low regularity}
Recall that in Section 4 we have constructed a solution $(u,\rho)$ of (\ref{eq:cNS}) satisfying
\ben\label{eq:properties}
\begin{split}
&\f {c_0} 2 \le \rho(t,x)\le 2c_0^{-1},\quad
\|a\|_{\widetilde{L}^\infty_T(\dot{B}^{
\frac{3}{p}}_{p,1})}+\|u\|_{\widetilde{L}^\infty_T(\dot B^{\frac
3p-1}_{p,1})} \le C_0E_0,\\
&\|a\|_{\widetilde{L}^\infty_T(\dot B^{
\frac{3}{p}}_{p,1}(\om))}\le \eta_1,
\quad \|u\|_{\widetilde{L}^1_T(\dot B^{ \frac{3}{p}+1}_{p,1})}
+\|u\|_{\widetilde{L}^2_T(\dot B^{ \frac{3}{p}}_{p,1})}\le \eta_2,
\end{split}
\een
with $a(t,x)=(\rho(t,x)-\overline{\rho})/\overline{\rho},\,E_0=\|a_0\|_{\dot B^{ \frac{3}{p}}_{p,1}}+\|u_0\|_{\dot B^{
\frac{3}{p}-1}_{p,1}}$.

We will show that this solution allows to propagate the regularity of the initial data in Sobolev space
with low regularity. This will ensure that one can use Hoff's theory for small energy initial data, since
the norm of the initial data is small in Sobolev space with low regularity,
but it is large in Sobolev space with high regularity.

\begin{Proposition}\label{prop:low}
Let $p\in (3,6)$ and $1-\f3 p<\delta <\f3p$. Assume that $(\rho,u)$ is a solution of (\ref{eq:cNS})
satisfying (\ref{eq:properties}). If $(a_0,u_0)\in \dot H^{1-\delta}\times \dot H^{-\delta}$, then
\beno
&&\|a\|_{\widetilde{L}^\infty_T(\dot B^{1-\delta}_{2,2})}
\leq C\big(\|a_0\|_{\dot H^{1-\delta}}+\|u_0\|_{\dot H^{-\delta}}\big),\\
&&\|u\|_{\widetilde{L}^\infty_T(\dot B^{-\delta}_{2,2})}\le C(1+E_0)\bigl(\|u_0\|_{\dot
H^{-\delta}}+T\|a_0\|_{\dot H^{1-\delta}}\bigr),\\
&&\|u\|_{\widetilde{L}^1_T(\dot B^{2-\delta}_{2,2})}+\|u\|_{\widetilde{L}^2_T(\dot B^{1-\delta}_{2,2})}
\le C\bigl(\|u_0\|_{\dot
H^{-\delta}}+T\|a_0\|_{\dot H^{1-\delta}}\bigr).
\eeno
Here the constant $C$ depends only on $\lambda,\mu, p,\delta, c_0$.
\end{Proposition}

\no{\bf Proof.\,} Set $\beta=1-\delta$. Without loss of generality, we can assume that
$\eta_1,\eta_2, T$ small enough such that
\ben
C\eta_1\le \f 18,\quad C\eta_2\le \f 18, \quad CT\le \f 18,
\een
where $C$ is the constant appearing in the following estimates and depending on
$\lambda,\mu, p,\delta, c_0$.

{\bf Step 1}. Estimate for the transport equation

Apply the operator $\Delta_j$ to the first equation of
(\ref{eq:linearcNS}) to obtain
\beno
\p_t\Delta_ja+u\cdot\na
\Delta_ja=\Delta_jF+[u,\Delta_j]\cdot \na a,
\eeno
with $F=-(1+a)\dv\,u.$ Making $L^2$ energy estimate, we get
\beno
\f12\f{d}{dt}\|\Delta_ja\|_2^2-\f12\int_{\R^3}|\Delta_ja|^{2}\dv
udx\le \bigl(\|\Delta_jF\|_2+\|[u,\Delta_j]\cdot \na
a\|_2\bigr)\|\Delta_ja\|_2,
\eeno
this gives
\beno
\|\Delta_ja(t)\|_2\le
\|\Delta_ja_0\|_2+\int_0^t\bigl(\|\Delta_jF\|_2+ \|[u,\Delta_j]\cdot
\na a\|_2+\f12\|\dv u\|_\infty\|\Delta_ja\|_2\bigr)d\tau,
\eeno
from which, we deduce that
\ben\label{a2}
\|a\|_{\widetilde{L}^\infty_t(\dot B^\beta_{2,2})}&\le& \|a_0\|_{\dot
B^\beta_{2,2}}+C\int^t_0\|\dv u\|_{\dot{B}^{\frac{3}{p}}_{p,1}}
\|a\|_{\tilde{L}^\infty_\tau(\dot{B}^\beta_{2,2})}d\tau\nonumber\\
&&+\|F\|_{\widetilde{L}^1_t(\dot{B}^\beta_{2,2})}+\big\|2^{j\beta}\|[u,\Delta_j]\cdot
\na a\|_{L_t^1(L^2)}\big\|_{\ell^2}.
\een
We get by Lemma \ref{Lem:mu2} and $\delta>1-\f 3p$ that
\beno
\|F\|_{\widetilde{L}^1_T(\dot B^\beta_{2,2})} \leq
C\|u\|_{\widetilde{L}^1_T(\dot
B^{\beta+1}_{2,2})}\big(1+\|a\|_{\widetilde{L}^\infty_T(\dot B^{\f3 p}_{p,1})}\big).
\eeno
and by Lemma \ref{Lem:com-Hs1},
\beno
&&\big\|2^{j\beta}\|[u,\Delta_j]\cdot \na a\|_{L_T^1(L^2)}\big\|_{\ell^2}
\le C\|u\|_{\tilde{L}_T^1(\dot
B^{\frac{3}{p}+1}_{p,1})}\|a\|_{\tilde{L}_T^\infty(\dot B^\beta_{2,2})}.
\eeno
Plugging them into (\ref{a2}) and using (\ref{eq:properties}),
we infer from Gronwall's inequality that
\ben\label{eq:a-low}
\|a\|_{\widetilde{L}^\infty_T(\dot B^\beta_{2,2})}
\leq C\big(\|a_0\|_{\dot
H^{1-\delta}}+\|u\|_{\widetilde{L}^1_T(\dot B^{\beta+1}_{2,2})}\big).
\een

{\bf Step 2.} Estimate for the momentum equation

Apply the operator $\Delta_j$ to the second equation of
(\ref{eq:linearcNS}) to obtain
\begin{eqnarray*}
\p_t\Delta_j u-\dv(\overline{\mu}\na
\Delta_ju)-\na(\overline{\lambda}\dv\,\Delta_ju)=\Delta_jG+\dv([\Delta_j,\overline{\mu}]\na
u) +\na([\overline{\lambda},\Delta_j]\dv{u}),
\end{eqnarray*}
and recall that
\beno
G=-u\cdot\na u-\f{\overline{\rho}P'(\rho)}
{\rho}\na a+\f {\mu} {\rho^2}\na \rho\cdot\na u+\f
{\lambda}{\rho^2}\na \rho\dv{u}.
\eeno
Making $L^2$ energy estimate(see the proof of Proposition \ref{Prop:momenequ}), we obtain
\begin{eqnarray}
\|u\|_{\widetilde{L}^\infty_T(\dot B^{\beta-1}_{2,2})}&\le&
C\bigl(\|u_0\|_{\dot B^{\beta-1}_{2,2}}
+\|G\|_{\widetilde{L}^1_T(\dot B^{\beta-1}_{2,2})}\bigr)\nonumber\\
&&+C\big\|2^{j\beta}\bigl(\|[\Delta_j,\overline{\nu}]\na
u\|_{L^1_T(L^2)}+\|[\Delta_j,\overline{\mu}]\na
u\|_{L^1_T(L^2)}\bigr)\big\|_{\ell^2},\label{eq:unweight}
\end{eqnarray}
and
\begin{eqnarray}\label{u}
&&\|u\|_{\widetilde{L}^1_T(\dot B^{\beta+1}_{2,2})}+\|u\|_{\widetilde{L}^2_T(\dot B^{\beta}_{2,2})}\nonumber\\
&&\le C\bigl(\|u_0\|_{\dot B^{\beta-1}_{2,2}(\om)}+\|G\|_{\widetilde{L}^1_T(\dot B^{\beta-1}_{2,2}(\om))}\bigr)\nonumber\\
&&\quad+C\big\|2^{j\beta}w_j(T)\bigl(\|[\Delta_j,\overline{\nu}]\na
u\|_{L^1_T(L^2)}+\|[\Delta_j,\overline{\mu}]\na
u\|_{L^1_T(L^2)}\bigr)\big\|_{\ell^2}.
\end{eqnarray}

First of all, by Lemma \ref{Lem:com-Hs}($\delta>1-\f3p$) and Lemma \ref{Lem:non-est}, we get
\begin{eqnarray}\label{commutor2}
&&\big\|2^{j\beta}\bigl(\|[\Delta_j,\overline{\nu}]\na
u\|_{L^1_T(L^2)}+\|[\Delta_j,\overline{\mu}]\na
u\|_{L^1_T(L^2)}\bigr)\big\|_{\ell^2}\nonumber\\
&&\leq C\|a\|_{\widetilde{L}^\infty_T (\dot
B^{\f3 p}_{p,1})}\|u\|_{\widetilde{L}^1_T(\dot B^{\beta+1}_{2,2})},
\end{eqnarray}
and
\begin{eqnarray}\label{commutor1}
&&\big\|2^{j\beta}w_j(T)\bigl(\|[\Delta_j,\overline{\nu}]\na
u\|_{L^1_T(L^2)}+\|[\Delta_j,\overline{\mu}]\na
u\|_{L^1_T(L^2)}\bigr)\big\|_{\ell^2}\nonumber\\
&&\leq C\|a\|_{\widetilde{L}^\infty_T (\dot
B^{\f3 p}_{p,1}(\om))}\|u\|_{\widetilde{L}^1_T(\dot B^{\beta+1}_{2,2})}.
\end{eqnarray}
Applying Lemma \ref{Lem:mu2} and Lemma \ref{Lem:non-est} to $G$, we obtain
\begin{eqnarray}\label{G1}
\|G\|_{\widetilde{L}^1_T(\dot B^{\beta-1}_{2,2})} &\leq& C\bigl(
\|u\|_{\widetilde{L}^2_T(\dot
B^{\f3 p}_{p,1})}\|u\|_{\widetilde{L}^2_T(\dot B^{\beta}_{2,2})}
+\|a\|_{\widetilde{L}^\infty_T(\dot B^{\f3 p}_{p,1})}\|u\|_{\widetilde{L}^1_T(\dot B^{\beta+1}_{2,2})}\nonumber\\
&&\qquad+T\|a\|_{\widetilde{L}^\infty_T(\dot B^{\beta}_{2,2})}\bigr),
\end{eqnarray}
and
\begin{eqnarray}\label{G}
\|G\|_{\widetilde{L}^1_T(\dot B^{\beta-1}_{2,2}(\om))} &\leq& C\bigl(
\|u\|_{\widetilde{L}^2_T(\dot
B^{\f3 p}_{p,1})}\|u\|_{\widetilde{L}^2_T(\dot B^{\beta}_{2,2})}
+\|a\|_{\widetilde{L}^\infty_T(\dot B^{\f3 p}_{p,1}(\om))}\|u\|_{\widetilde{L}^1_T(\dot B^{\beta+1}_{2,2})}\nonumber\\
&&\qquad+T\|a\|_{\widetilde{L}^\infty_T(\dot B^{\beta}_{2,2})}\bigr).
\end{eqnarray}
Here we need $1-\f 3 p<\delta<\f 3p$.

Plugging (\ref{commutor2}) and (\ref{G1}) into (\ref{eq:unweight}) and using (\ref{eq:properties}),
we get
\ben\label{eq:u-infty}
&&\|u\|_{\widetilde{L}^\infty_T(\dot B^{\beta-1}_{2,2})}\le C\bigl(\|u_0\|_{\dot
H^{-\delta}}+T\|a\|_{\widetilde{L}^\infty_T(\dot
B^{\beta}_{2,2})}+E_0\|u\|_{\widetilde{L}^1_T(\dot B^{\beta+1}_{2,2})}\bigr),
\een
and plugging (\ref{commutor1}) and (\ref{G}) into (\ref{u}), we have
\begin{eqnarray}\label{eq:u-L1}
&&\|u\|_{\widetilde{L}^1_T(\dot B^{\beta+1}_{2,2})}+\|u\|_{\widetilde{L}^2_T(\dot B^{\beta}_{2,2})}
\le C\bigl(\|u_0\|_{\dot
H^{-\delta}}+T\|a\|_{\widetilde{L}^\infty_T(\dot
B^{\beta}_{2,2})}\bigr)
\end{eqnarray}
Then Proposition \ref{prop:low} follows by combining (\ref{eq:u-infty})-(\ref{eq:u-L1}) with
(\ref{eq:a-low}).\ef

\subsection{Propagation of high regularity}

We will show that the solution constructed in Section 4 can also propagate the regularity of the initial data in Sobolev space
with high regularity.

\begin{Proposition}\label{prop:high}
Assume that $(\rho,u)$ is a solution of (\ref{eq:cNS}) on $[0,T]$, which satisfies
\ben\label{eq:(a,u)-class}
\rho\ge c_0,\quad a\in \widetilde{L}^\infty(0,T;\dot B^{\f 3p}_{p,1}),\quad
u\in \widetilde{L}^\infty(0,T;\dot B^{\f 3p-1}_{p,1})\cap \widetilde{L}^1(0,T;\dot B^{\f 3p+1}_{p,1}).
\een
If $(a_0,u_0)\in H^s\times H^{s-1}$ for $s\ge 3$, then we have
\beno
a\in \widetilde{L}^\infty(0,T;\dot B^{s}_{2,2}),\quad
u\in \widetilde{L}^\infty(0,T;\dot B^{s-1}_{2,2})\cap \widetilde{L}^1(0,T;\dot B^{s+1}_{2,2}).
\eeno

\end{Proposition}

\no{\bf Proof.}\,Due to (\ref{eq:(a,u)-class}), we can divide the time interval $[0,T]$ into finite many small intervals
$[T_i, T_{i+1}]$ with $i=0,\cdots, N$ such that
\beno
&&\|a\|_{\widetilde{L}^\infty(T_i,T_{i+1};\dot B^{\f 3p}_{p,1}(\om^i))}\le \epsilon,\\
&&\|u\|_{\widetilde{L}^1(T_i,T_{i+1};\dot B^{\f 3p+1}_{p,1})}+\|u\|_{\widetilde{L}^2(T_i,T_{i+1};\dot B^{\f 3p}_{p,1})}\le \epsilon,
\eeno
for some $\epsilon$ small enough, see Remark \ref{rem:interval}. Here the weight on each $[T_i,T_{i+1}]$ is given by
$$
\om_k^i=\sum_{\ell\ge k}2^{k-\ell}e_\ell^i,\quad e_\ell^i=(1-e^{-c2^{2\ell}(T_{i+1}-T_i)})^\f12.
$$

First of all, we also have(see \ref{a2})
\ben\label{eq:a-high}
\|a\|_{\widetilde{L}^\infty_t(\dot B^s_{2,2})}&\le& \|a_0\|_{\dot
B^s_{2,2}}+C\int^t_0\|\dv u\|_{\dot{B}^{\frac{3}{p}}_{p,1}}
\|a\|_{\tilde{L}^\infty_\tau(\dot{B}^s_{2,2})}d\tau\nonumber\\
&&+\|F\|_{\widetilde{L}^1_t(\dot{B}^s_{2,2})}+\big\|2^{js}\|[u,\Delta_j]\cdot
\na a\|_{L_t^1(L^2)}\big\|_{\ell^2}.
\een
By Lemma \ref{Lem:prod-H}, we have
\beno
\|F\|_{\widetilde{L}^1_t(\dot B^s_{2,2})}
\leq C\big(\|u\|_{\widetilde{L}^1_t(\dot B^{s+1}_{2,2})}+\|a\|_{\widetilde{L}^{\infty}_t(\dot
B^{\f3 p}_{p,1})}\|u\|_{\widetilde{L}^{1}_t(\dot B^{s+1}_{2,2})}+\|a\|_{\widetilde{L}^{\infty}_t(\dot
B^{s}_{2,2})}\|u\|_{\widetilde{L}^{1}_t(\dot B^{\f3p+1}_{p,1})}\big),
\eeno
and by Lemma \ref{Lem:com-H2},
\beno
\big\|2^{js}\|[u,\Delta_j]\cdot \na a\|_{L_t^1(L^2)}\big\|_{\ell^2} \le {C}\big(
\|u\|_{\tilde{L}_t^1(\dot{B}_{p,1}^{\frac{3}{p}+1})}\|a\|_{\tilde{L}_t^\infty(\dot{B}_{2,2}^s)}
+\|u\|_{\tilde{L}_t^1(\dot{B}_{2,2}^{s+1})}\|a\|_{\tilde{L}_t^\infty(\dot{B}_{p,1}^{\frac{3}{p}})}\big),
\eeno
from which and (\ref{eq:a-high}), we get by Gronwall's inequality that for $t\in [0,T_1]$,
\ben\label{eq:a-high-energy}
\|a\|_{\widetilde{L}^\infty_t(\dot B^s_{2,2})}\leq C\big(\|a_0\|_{\dot{B}^s_{2,2}}+
\|u\|_{\widetilde{L}^1_t(\dot B^{s+1}_{2,2})}\big).
\een

Making $H^{s-1}$ energy estimate for the momentum equation, we get(see \ref{eq:unweight} and (\ref{u}))
\begin{eqnarray}\label{eq:u-high1}
\|u\|_{\widetilde{L}^\infty_t(\dot B^{s-1}_{2,2})}&\le&
C\bigl(\|u_0\|_{\dot B^{s-1}_{2,2}}
+\|G\|_{\widetilde{L}^1_t(\dot B^{s-1}_{2,2})}\bigr)\nonumber\\
&&+C\big\|2^{js}\bigl(\|[\Delta_j,\overline{\nu}]\na
u\|_{L^1_t(L^2)}+\|[\Delta_j,\overline{\mu}]\na
u\|_{L^1_t(L^2)}\bigr)\big\|_{\ell^2},
\end{eqnarray}
and
\begin{eqnarray}\label{eq:u-high2}
&&\|u\|_{\widetilde{L}^1_t(\dot B^{s+1}_{2,2})}+\|u\|_{\widetilde{L}^2_t(\dot B^{s}_{2,2})}\nonumber\\
&&\le C\bigl(\|u_0\|_{\dot B^{s-1}_{2,2}}+\|G\|_{\widetilde{L}^1_t(\dot B^{s-1}_{2,2}(\om))}\bigr)\nonumber\\
&&\quad+C\big\|2^{js}w_j(t)\bigl(\|[\Delta_j,\overline{\nu}]\na
u\|_{L^1_t(L^2)}+\|[\Delta_j,\overline{\mu}]\na
u\|_{L^1_t(L^2)}\bigr)\big\|_{\ell^2}.
\end{eqnarray}
By Lemma \ref{Lem:prod-H} and Lemma \ref{Lem:non-est}, we have
\beno
&&\|G\|_{\widetilde{L}^1_t(\dot B^{s-1}_{2,2}(\om))}
\leq C\Bigl( \|u\|_{\widetilde{L}^2_t(\dot B^{\f3p}_{p,1})}\|u\|_{\widetilde{L}^2_t(\dot B^{s}_{2,2})}
+\|u\|_{\widetilde{L}^2_t(\dot B^{\f3p}_{p,1})}\|u\|_{\widetilde{L}^2_T(\dot B^{s}_{2,2})}\nonumber\\
&&\quad+\|a\|_{\widetilde{L}^\infty_t(\dot B^{\f3p}_{p,1}(\om))}\|u\|_{\widetilde{L}^1_t(\dot B^{s+1}_{2,2})}
+\|a\|_{\widetilde{L}^\infty_t(\dot B^{s}_{2,2})}\|u\|_{\widetilde{L}^1_t(\dot B^{\f3p+1}_{p,1})}
+t\|a\|_{\widetilde{L}^\infty_t(\dot B^{s}_{2,2})}\Bigr),
\eeno
and Lemma \ref{Lem:com-H1} and Lemma \ref{Lem:non-est},
\beno
&&\big\|2^{js}w_j(t)\bigl(\|[\Delta_j,\overline{\nu}]\na
u\|_{L^1_t(L^2)}+\|[\Delta_j,\overline{\mu}]\na
u\|_{L^1_t(L^2)}\bigr)\big\|_{\ell^2}\nonumber\\
&&\quad\leq C\Big(\|a\|_{\tilde{L}_t^\infty(\dot{B}_{p,1}^{\frac{3}{p}}(\om))}\|u\|_{\tilde{L}_t^1(\dot{B}_{2,2}^{s+1})}
+\|a\|_{\tilde{L}_t^{\infty}(\dot{B}_{2,2}^s)}\|u\|_{\tilde{L}_t^1(\dot{B}_{p,1}^{\frac{3}{p}+1})}\Big),
\eeno
which along with (\ref{eq:u-high2}) gives that
\beno
\|u\|_{\widetilde{L}^1_t(\dot B^{s+1}_{2,2})}+\|u\|_{\widetilde{L}^2_t(\dot B^{s}_{2,2})}
\le C\bigl(\|u_0\|_{\dot H^{s-1}}+\|a\|_{\widetilde{L}^\infty_t(\dot B^{s}_{2,2})}(\|u\|_{\widetilde{L}^1_t(\dot B^{\f3p+1}_{p,1})}+t)\big),
\eeno
from which and (\ref{eq:a-high-energy}), it follows that for any $t\in [0,T_1]$,
\ben\label{eq:(a,u)-high}
\|a\|_{\widetilde{L}^\infty_t(\dot B^s_{2,2})}+\|u\|_{\widetilde{L}^1_t(\dot B^{s+1}_{2,2})}+\|u\|_{\widetilde{L}^2_t(\dot B^{s}_{2,2})}
\le C\big(\|a_0\|_{\dot H^s}+\|u_0\|_{\dot H^{s-1}}\big).
\een
On the other hand, by Lemma \ref{Lem:prod-H} and Lemma \ref{Lem:non-est} again, we have
\beno
&&\|G\|_{\widetilde{L}^1_t(\dot B^{s-1}_{2,2})}
\leq C\Bigl( \|u\|_{\widetilde{L}^2_t(\dot B^{\f3p}_{p,1})}\|u\|_{\widetilde{L}^2_t(\dot B^{s}_{2,2})}
+\|u\|_{\widetilde{L}^2_t(\dot B^{\f3p}_{p,1})}\|u\|_{\widetilde{L}^2_T(\dot B^{s}_{2,2})}\nonumber\\
&&\quad+\|a\|_{\widetilde{L}^\infty_t(\dot B^{\f3p}_{p,1})}\|u\|_{\widetilde{L}^1_t(\dot B^{s+1}_{2,2})}
+\|a\|_{\widetilde{L}^\infty_t(\dot B^{s}_{2,2})}\|u\|_{\widetilde{L}^1_t(\dot B^{\f3p+1}_{p,1})}
+t\|a\|_{\widetilde{L}^\infty_t(\dot B^{s}_{2,2})}\Bigr),
\eeno
and by Lemma \ref{Lem:com-H1} and Lemma \ref{Lem:non-est},
\beno
&&\big\|2^{js}\bigl(\|[\Delta_j,\overline{\nu}]\na
u\|_{L^1_t(L^2)}+\|[\Delta_j,\overline{\mu}]\na
u\|_{L^1_t(L^2)}\bigr)\big\|_{\ell^2}\nonumber\\
&&\quad\leq C\Big(\|a\|_{\tilde{L}_t^\infty(\dot{B}_{p,1}^{\frac{3}{p}})}\|u\|_{\tilde{L}_t^1(\dot{B}_{2,2}^{s+1})}
+\|a\|_{\tilde{L}_t^{\infty}(\dot{B}_{2,2}^s)}\|u\|_{\tilde{L}_t^1(\dot{B}_{p,1}^{\frac{3}{p}+1})}\Big),
\eeno
which along with (\ref{eq:u-high1}) and (\ref{eq:a-high-energy}) implies that for any $t\in [0,T_1]$,
\beno
\|u\|_{\widetilde{L}^\infty_t(\dot B^{s-1}_{2,2})}
\le C\big(1+\|a\|_{\tilde{L}_t^\infty(\dot{B}_{p,1}^{\frac{3}{p}})}\big)\big(\|a_0\|_{\dot H^s}+\|u_0\|_{\dot H^{s-1}}\big).
\eeno
Then the proposition follows by repeating the above process in each interval $[T_i,T_{i+1}]$.\ef

\begin{Remark}\label{rem:interval}
Thanks to $a\in \widetilde{L}^\infty(0,T;\dot B^{\f 3p}_{p,1})$, one can choose $N_1$ big enough so that
\beno
\sum_{|j|>N_1}2^{j\f 3 p}\|\Delta_j a(t)\|_{L^\infty(0,T;L^p)}\le \f \epsilon 4.
\eeno
Fix $N_1$, one can divide $[0,T]$ into finite many small intervals
$[T_i, T_{i+1}]$ with $i=0,\cdots, N$ such that for $i=0,\cdots, N$
\beno
\sum_{|j|\le N_1}2^{j\f 3 p}\om_j^i\|\Delta_j a(t)\|_{L^\infty(T_i,T_{i+1};L^p)}\le \f \epsilon 2.
\eeno
Then we conclude that
\beno
&&\|a\|_{\widetilde{L}^\infty(T_i,T_{i+1};\dot B^{\f 3p}_{p,1}(\om^i))}\\
&&\le 2\sum_{|j|>N_1}2^{j\f 3 p}\|\Delta_j a(t)\|_{L^\infty(T_i,T_{i+1};L^p)}
+\sum_{|j|\le N_1}2^{j\f 3 p}\om_j^i\|\Delta_j a(t)\|_{L^\infty(T_i,T_{i+1};L^p)}\le \epsilon.
\eeno

\end{Remark}

%
%
%
%

\subsection{Some technical lemmas}

\begin{Lemma}\label{lem:mu}
Let $p\ge 2, 1\leq q_1,q_2,q\leq\infty$ with $\f 1{q_1}+\f 1{q_2}=\f{1}{q}$. Then there hold

\no (a)\, if $s_2< \frac{3}{p}$, we have
\beno \|T_gf\|_{\widetilde{L}^q_T(\dot
B^{s_1+s_2- \frac{3}{p}}_{2,2}(\omega))} \le C\|f\|_{\widetilde{L}^{q_1}_T(\dot
B^{s_1}_{p,1}(\omega))}\|g\|_{\widetilde{L}^{q_2}_T(\dot B^{s_2}_{2,2})};
\eeno

\no (b)\, if $s_1\le \frac{3}{p}-1$, we have
\beno
\|T_fg\|_{\widetilde{L}^q_T(\dot
B^{s_1+s_2- \frac{3}{p}}_{2,2}(\omega))}\le C\|f\|_{\widetilde{L}^{q_1}_T(\dot
B^{s_1}_{p,1}(\omega))}\|g\|_{\widetilde{L}^{q_2}_T(\dot B^{s_2}_{2,2})};
\eeno

\no (c)\, if $s_1+s_2>0$,
we have
\beno
\|R(f,g)\|_{\widetilde{L}^q_T(\dot
B^{s_1+s_2- \frac{3}{p}}_{2,2}(\omega))} \le
C\|f\|_{\widetilde{L}^{q_1}_T(\dot
B^{s_1}_{p,1}(\omega))}\|g\|_{\widetilde{L}^{q_2}_T(\dot B^{s_2}_{2,2})}.
\eeno
\end{Lemma}
\no{\bf Proof.}\, We first prove (a). Due to (\ref{orth}), we have \beno
\Delta_j(T_gf)=\sum_{|j'-j|\le 4}\Delta_j(S_{j'-1}g\Delta_{j'}f),
\eeno then we get by Lemma \ref{Lem:Bernstein} and
(\ref{weightprop}) that
\beno
\|T_gf\|_{\widetilde{L}^q_T(\dot
B^{s_1+s_2- \frac{3}{p}}_{2,2}(\omega))}&=&\big\|2^{j(s_1+s_2- \frac{3}{p})}\om_j(T)\|\Delta_j(T_gf)\|_{L^q_T(L^2)}\big\|_{\ell^2}\\
&\le& C\big\|2^{j(s_1+s_2- \frac{3}{p})}\om_{j}(T)\|S_{j-1}g\|_{L^{q_2}_T(L^{\frac{2p}{p-2}})}\|\Delta_{j}f\|_{L^{q_1}_T(L^p)}\big\|_{\ell^2}\\
&\le& C\|f\|_{\widetilde{L}^{q_1}_T(\dot
B^{s_1}_{p,1}(\omega))}\|g\|_{\widetilde{L}^{q_2}_T(\dot B^{s_2}_{2,2})},
\eeno
where we used in the last inequality
\beno
\|S_{j-1}g\|_{L^{q_2}_T(L^{\frac{2p}{p-2}})} \le C\sum_{\ell\le j-2}2^{\ell
\frac{3}{p}}\|\Delta_{\ell}g\|_{L^{q_2}_T(L^{2})}
\le C2^{j(-s_2+\frac{3}{p})}\|g\|_{\widetilde{L}^{q_2}_T(\dot B^{s_2}_{2,2})}.
\eeno

We next prove (b). Similarly,
we have
\beno
\|T_fg\|_{\widetilde{L}^q_T(\dot
B^{{s_1+s_2- \frac{3}{p}}}_{2,2}(\omega))}&=&\big\|2^{j({s_1+s_2- \frac{3}{p}})}\om_j(T)\|\Delta_j(T_fg)\|_{L^q_T(L^2)}\big\|_{\ell^2}\\
&\le& C\big\|2^{j({s_1+s_2- \frac{3}{p}})}\om_{j}(T)\|S_{j-1}f\|_{L^{q_1}_T(L^\infty)}\|\Delta_{j}g\|_{L^{q_2}_T(L^2)}\big\|_{\ell^2},
\eeno
and by Lemma \ref{Lem:Bernstein} and (\ref{weightprop}), we
have
\beno
\om_{j}(T)\|S_{j-1}f\|_{L^{q_1}_T(L^\infty)}&\le& C2^{j}\sum_{\ell\le j-2}2^{\ell ( \frac{3}{p}-1)}\om_{\ell}(T)\|\Delta_{\ell}f\|_{L^{q_1}_T(L^p)}\\
&\le&C2^{j(\frac{3}{p}-s_1)}\sum_{\ell\le{j-2}}2^{\ell{s_1}}\om_{\ell}(T)\|\Delta_{\ell}f\|_{L^{q_1}_T(L^p)}\\
&\le&C2^{j(\frac{3}{p}-s_1)}\|f\|_{\widetilde{L}^{q_1}_T(\dot
B^{s_1}_{p,1}(\omega))},
\eeno
 which lead to (b). Finally, let us prove (c). Noticing that
 \beno
\Delta_j(R(f,g))=\sum_{j'\ge
j-3}\Delta_j(\Delta_{j'}f\widetilde{\Delta}_{j'}g),
\eeno
then we get by Lemma \ref{Lem:Bernstein} that
\beno
&&\|R(f,g)\|_{\widetilde{L}^q_T(\dot
B^{{s_1+s_2-\frac{3}{p}}}_{2,2}(\omega))}\\
&&\le C\big\|\sum_{j'\ge j-3}2^{j(s_1+s_2)}\om_j(T)\|\Delta_{j'}f\|_{L^{q_1}_T(L^p)}\|\widetilde{\Delta}_{j'}g\|_{L^{q_2}_T(L^2)}\big\|_{\ell^2}\\
&&\le C\big\|\sum_{j'\ge j-3}\sum_{\ell\ge j}2^{j-\ell}e_\ell(T)2^{j(s_1+s_2)}\|\Delta_{j'}f\|_{L^{q_1}_T(L^p)}\|\widetilde{\Delta}_{j'}g\|_{L^{q_2}_T(L^2)}\big\|_{\ell^2}\\
&&= \big\|\sum_{j'\ge j-3}\sum_{\ell\ge j,j'}\cdots\big\|_{\ell^2}
+\big\|\sum_{j'\ge j-3}\sum_{j'\ge \ell\ge j}\cdots\big\|_{\ell^2}\\
&&=I+II.
 \eeno
In the following, we denote $\{c_j\}$ by a sequence in $\ell^2$.
Noting that
\beno \sum_{\ell\ge
j,j'}2^{j-\ell}e_\ell(T)\le 2^{j-j'}\sum_{\ell\ge
j'}2^{j'-\ell}e_\ell(T)=2^{j-j'}w_{j'}(T),
\eeno
we infer that
\beno
I&\le& C\Big\|\sum_{j'\ge j-3}\om_{j'}(T)2^{j(s_1+s_2)}2^{j-j'}\|\Delta_{j'}f\|_{L^{q_1}_T(L^p)}\|\widetilde{\Delta}_{j'}g\|_{L^{q_2}_T(L^2)}\Big\|_{\ell^2}\\
&\le& C\Big\|\sum_{j'\ge j-3}\om_{j'}(T)2^{(j-j')(s_1+s_2+1)}2^{j's_1}c_{j'}\|\Delta_{j'}f\|_{L^{q_1}_T(L^p)}\Big\|_{\ell^2}\|g\|_{\widetilde{L}^{q_2}_T(\dot B^{s_2}_{2,2})}\\
&\le& C\|f\|_{\widetilde{L}^{q_1}_T(\dot
B^{s_1}_{p,1}(\omega))}\|g\|_{\widetilde{L}^{q_2}_T(\dot B^{s_2}_{2,2})}.
\eeno
and thanks to $e_\ell(T)\le e_{j'}(T)\le \om_{j'}(T)$ for $\ell\le j'$, we get
\beno
II&\le&C \Big\| \sum_{j'\ge j-3}\om_{j'}(T)2^{j(s_1+s_2)}\|\Delta_{j'}f\|_{L^{q_1}_T(L^p)}\|\widetilde{\Delta}_{j'}g\|_{L^{q_2}_T(L^2)}\Big\|_{\ell^2}\\
&\le& C \Big\| \sum_{j'\ge j-3}\om_{j'}(T)2^{(j-j')(s_1+s_2)}2^{j's_1}c_{j'}\|\Delta_{j'}f\|_{p}\Big\|_{\ell^2}\|g\|_{\widetilde{L}^{q_2}_T(\dot B^{s_2}_{2,2})}\\
&\le& C\|f\|_{\widetilde{L}^{q_1}_T(\dot
B^{s_1}_{p,1}(\omega))}\|g\|_{\widetilde{L}^{q_2}_T(\dot B^{s_2}_{2,2})},
\eeno
 This proves (c). \ef

From Lemma \ref{lem:mu} and its proof, it is easy to see that
\begin{Lemma}\label{Lem:mu2}
Let $p\ge 2,\,s_1\le \frac{3}{p},\,s_2<\f 3p,\, s_1+s_2>0$, and $1\le q,q_1,q_2\le \infty$ with $\f 1{q_1}+\f1{q_2}=\f1q$.
Then we have
\beno
\|fg\|_{\widetilde{L}^{q}_T(\dot B^{s_1+s_2-
\frac{3}{p}}_{2,2})} \le C\|f\|_{\widetilde{L}^{q_1}_T(\dot
B^{s_1}_{p,1})}\|g\|_{\widetilde{L}^{q_2}_T(\dot B^{s_2}_{2,2})},
\eeno
and if $s_1\le \f 3 p-1$, then
\beno \|fg\|_{\widetilde{L}^{q}_T(\dot B^{s_1+s_2-
\frac{3}{p}}_{2,2}(\om))} \le C\|f\|_{\widetilde{L}^{q_1}_T(\dot
B^{s_1}_{p,1}(\om))}\|g\|_{\widetilde{L}^{q_2}_T(\dot B^{s_2}_{2,2})}.
\eeno
\end{Lemma}

\begin{Lemma}\label{Lem:com-Hs}
Let $p\ge 2$ and $s\in(-\frac{3}{p},\frac{3}{p})$. Then there holds
\beno
&&\big\|2^{js}\|[\Delta_j,f]\na g\|_{L_T^1(L^2)}\big\|_{\ell^2}\le{C}
\|f\|_{\tilde{L}_T^\infty(\dot{B}_{p,1}^{\frac{3}{p}})}\|g\|_{\tilde{L}_T^1(\dot{B}_{2,2}^{s+1})},\\
&&\big\|2^{js}\om_j(T)\|[\Delta_j,f]\na g\|_{L_T^1(L^2)}\big\|_{\ell^2}\le{C}
\|f\|_{\tilde{L}_T^\infty(\dot{B}_{p,1}^{\frac{3}{p}}(\om))}\|g\|_{\tilde{L}_T^1(\dot{B}_{2,2}^{s+1})}.
\eeno
\end{Lemma}

\no{\bf Proof.}\,The proof is almost the same as Lemma \ref{Lem:com-est}. We write
\beno
[f,\Delta_j]\na g=[T_f,\Delta_j]\na g+T_{\Delta_j\na g}f+R(f,\Delta_j\na g)-\Delta_j(T_{\na g}f)-\Delta_jR(f,\na g).
\eeno
It follows from Lemma \ref{lem:mu} (a), (c) that
\beno
&&\big\|2^{js}\om_j(T)\|\Delta_j(T_{\na g}f)\|_{L_T^1(L^2)}\big\|_{\ell^2}
\le C\|f\|_{\tilde{L}_T^\infty(\dot{B}_{p,1}^{\frac{3}{p}}(\om))}\|g\|_{\tilde{L}_T^1(\dot{B}_{2,2}^{s+1})},\\
&&\big\|2^{js}\om_j(T)\|\Delta_jR(f,\na g)\|_{L_T^1(L^2)}\big\|_{\ell^2}
\le C\|f\|_{\tilde{L}_T^\infty(\dot{B}_{p,1}^{\frac{3}{p}}(\om))}\|g\|_{\tilde{L}_T^1(\dot{B}_{2,2}^{s+1})}.
\eeno
Recalling (\ref{eq:para}) and by Lemma \ref{Lem:Bernstein} and (\ref{weightprop}), we get
\beno
&&\big\|\om_j(T)2^{js}\|T_{\Delta_j\na g}'f\|_{L^1_T(L^2)}\big\|_{\ell^2}\\
&&\le C\big\|\om_j(T)2^{j(s+1)}\|\Delta_jg\|_{L^1_T(L^2)}\sum_{j'\ge j-2}\|\Delta_{j'}f\|_{L^\infty_T(L^\infty)}\big\|_{\ell^2}\\
&&\le C\big\|2^{j(s+1)}\|\Delta_jg\|_{L^1_T(L^2)}\sum_{j'\ge j-2}\om_{j'}(T)2^{j'\frac{3}{p}}\|\Delta_{j'}f\|_{L^\infty_T(L^p)}\big\|_{\ell^2}\\
&&\le C\|f\|_{\widetilde{L}^\infty_T(\dot B^{
\frac{3}{p}}_{p,1}(\om))}\|g\|_{\widetilde{L}^1_T(\dot
B^{s+1}_{2,2})}.
\eeno
Note that
\beno [T_{f},
\Delta_j]\pa_k g
&=&\sum_{|j'-j|\le4}2^{4j}\int_{\R^3}\int_0^1y\cdot\na
S_{j'-1}f(x-\tau
y)d\tau\pa_kh(2^jy)\Delta_{j'}g(x-y)dy\nonumber\\&&\qquad\quad+
2^{3j}\int_{\R^3}h(2^j(x-y))\pa_kS_{j'-1}f(y)\Delta_{j'}g(y)dy ,
\eeno
from which and Young's inequality, we infer that
\beno
\big\|\om_j(T)2^{js}\|[T_{f}, \Delta_j]\na
g\|_{L^1_T(L^2)}\big\|_{\ell^2}\le C\|f\|_{\widetilde{L}^\infty_T(\dot B^{
\frac{3}{p}}_{p,1}(\om))}\|g\|_{\widetilde{L}^1_T(\dot
B^{s+1}_{2,2})}.
\eeno

Summing up the above estimates yield the second inequality, the first one is similar.\ef

From the proof of Lemma \ref{Lem:com-Hs}, it is easy to see that

\begin{Lemma}\label{Lem:com-Hs1}
Let $p\ge 2$ and $s\in(-\frac{3}{p},\frac{3}{p}+1)$. Then there holds
\beno
\big\|2^{js}\|[\Delta_j,f]\nabla{g}\|_{L_T^1(L^2)}\big\|_{\ell^2}\le{C}
\|f\|_{\tilde{L}_T^\infty(\dot{B}_{p,1}^{\frac{3}{p}+1})}\|g\|_{\tilde{L}_T^1(\dot{B}_{2,2}^s)}.
\eeno
\end{Lemma}

The following lemmas are used in the proof of propagation of high regularity. Since the proof is very similar
to Lemma \ref{Lem:mu2}-Lemma \ref{Lem:com-Hs1}, we state them without proof.

\begin{Lemma}\label{Lem:prod-H}
Let $p\ge 2, s_1, s_1'\le \f 3p, s_1+s_2=s_1'+s_2'>0$, and $1\le q,q_1,q_2\le \infty$ with $\f 1{q_1}+\f1{q_2}=\f1q$.
Then there holds
\beno
\|fg\|_{\widetilde{L}^{q}_T(\dot B^{s_1+s_2-\f 3p}_{2,2})} \le C\big(\|f\|_{\widetilde{L}^{q_1}_T(\dot
B^{s_1}_{p,1})}\|g\|_{\widetilde{L}^{q_2}_T(\dot B^{s_2}_{2,2})}+\|f\|_{\widetilde{L}^{q_1}_T(\dot
B^{s_2'}_{2,2})}\|g\|_{\widetilde{L}^{q_2}_T(\dot B^{s_1'}_{p,1})}\big),
\eeno
and if $s_1, s_1'\le \f 3p-1$, then
\beno
\|fg\|_{\widetilde{L}^{q}_T(\dot B^{s_1+s_2-\f 3p}_{2,2}(\om))} \le C\big(\|f\|_{\widetilde{L}^{q_1}_T(\dot
B^{s_1}_{p,1}(\om))}\|g\|_{\widetilde{L}^{q_2}_T(\dot B^{s_2}_{2,2})}+\|f\|_{\widetilde{L}^{q_1}_T(\dot
B^{s_2'}_{2,2})}\|g\|_{\widetilde{L}^{q_2}_T(\dot B^{s_1'}_{p,1}(\om))}\big).
\eeno
\end{Lemma}

\begin{Lemma}\label{Lem:com-H1}
Let $p\ge 2$ and $s>-\frac{3}{p}$. Then there holds
\beno
&&\big\|2^{js}\|[\Delta_j,f]\na g\|_{L_T^1(L^2)}\big\|_{\ell^2}\le{C}\big(
\|f\|_{\tilde{L}_T^\infty(\dot{B}_{p,1}^{\frac{3}{p}}(\om))}\|g\|_{\tilde{L}_T^1(\dot{B}_{2,2}^{s+1})}
+\|f\|_{\tilde{L}_T^{\infty}(\dot{B}_{2,2}^s(\om))}\|g\|_{\tilde{L}_T^1(\dot{B}_{p,1}^{\frac{3}{p}})}\big),\\
&&\big\|2^{js}\om_j(T)\|[\Delta_j,f]\na g\|_{L_T^1(L^2)}\big\|_{\ell^2}\le{C}\big(
\|f\|_{\tilde{L}_T^\infty(\dot{B}_{p,1}^{\frac{3}{p}}(\om))}\|g\|_{\tilde{L}_T^1(\dot{B}_{2,2}^{s+1})}
+\|f\|_{\tilde{L}_T^{\infty}(\dot{B}_{2,2}^s(\om))}\|g\|_{\tilde{L}_T^1(\dot{B}_{p,1}^{\frac{3}{p}})}\big).
\eeno

\end{Lemma}

\begin{Lemma}\label{Lem:com-H2}
Let $p\ge 2$ and $s>-\frac{3}{p}$. Then there holds
\beno
\big\|2^{js}\|[\Delta_j,f]\nabla{g}\|_{L_T^1(L^2)}\big\|_{\ell^2}\le{C}\big(
\|f\|_{\tilde{L}_T^1(\dot{B}_{p,1}^{\frac{3}{p}+1})}\|g\|_{\tilde{L}_T^\infty(\dot{B}_{2,2}^s)}
+\|f\|_{\tilde{L}_T^{1}(\dot{B}_{2,2}^{s+1})}\|g\|_{\tilde{L}_T^\infty(\dot{B}_{p,1}^{\frac{3}{p}})}\big).
\eeno
\end{Lemma}

\setcounter{equation}{0}
\section{Hoff's energy method}

In a series of papers \cite{Hoff-JDE, Hoff-ARMA, Hoff-CPAM, Hoff-JMFM}, Hoff developed a method to construct global weak solution
for the discontinuous initial data with small energy. An important property of this solution is that  the density does not develop vacuum and can not concentrate
if the initial density is bounded below and above. In this section, we will adapt his method to our case.
Due to the stronger condition on the initial data, the restriction on the viscosity coefficient can be
removed and the obtained estimate is also better.

We set $\sigma(t)=\min(1,t)$ and define
\beno
&&A_1(T)=\sup_{t\in [0,T]}\|\nabla u(t)\|^2_{L^2}+\int_0^T\int_{\R^3}\rho|\dot u|^2dx dt,\\
&&A_2(T)=\sup_{t\in [0,T]}\sigma(t)\int_{\R^3}\rho|\dot u(t,x)|^2dx+\int_0^T\int_{\R^3}\sigma|\nabla \dot u|^2dxdt.
\eeno
Here and what follows, we denote
\beno
\dot f=f_t+u\cdot\na f.
\eeno

\begin{Theorem}\label{thm:hoff}
Let $(\rho,u)$ be a solution of (\ref{eq:cNS}) satisfying
\beno
\rho-\overline{\rho} \in C([0,T];H^2),\quad u\in C([0,T];H^2)\cap L^2(0,T;H^3).
\eeno
There exists a constant $\varepsilon_0$ depending only on $\lambda, \mu, c_0, \overline{\rho}$ such that
if the initial data $(\rho_0,u_0)$ satisfies
\beno
&&c_0\le \rho_0(x)\le c_0^{-1},\quad x\in \R^3,\\
&&\int_{\R^3}|\rho_0(x)-\overline{\rho}|^2dx+\|u_0\|_{H^1}^2=C_0\le \varepsilon_0,
\eeno
then we have
\beno
&& \f {c_0}2 \le \rho(t,x)\le 2c_0^{-1},\quad (t,x)\in [0,T]\times \R^3,\\
&&A_1(T)+ A_2(T)\le C_0^\f12.
\eeno
\end{Theorem}

\no{\bf Proof.}\,Due to the assumption, there exists a $0<T_0\le T$
such that the solution $(\rho,u)$ satisfies
\beno
&&\f {c_0} 2\le \rho(t,x)\le 2c_0^{-1},\quad (t,x)\in [0,T_0]\times \R^3,\\
&&A_1(T_0)+ A_2(T_0)\le C_0^\f12.
\eeno
Without loss of generality, we assume that $T_0$ is a maximal time so that the above inequalities hold.
In the following, we will give a refined estimate on $[0,T_0]$ for the solution. We denote $C$ by a constant
depending only on $\lambda,\mu, c_0,\overline{\rho}$.

{\bf Step 1.}\,$L^2$ energy estimate
\ben\label{eq:L2}
&&\int_{\R^3}\f 12\rho(t,x)|u(t,x)|^2+G(\rho(t,x))dx+\int_0^t\int_{\R^3}|\na u|^2dxdt\nonumber\\
&&\leq \int_{\R^3}\f 12\rho_0|u_0|^2+G(\rho_0)dx\leq CC_0,
\een
where $G(\rho)$ is the potential energy density defined by
\beno
G(\rho)=\rho\int_{\overline{\rho}}^{\rho}\f {P(s)-P(\overline{\rho})} {s^2}ds.
\eeno
Due to (\ref{ass:pressure}), we have
\ben\label{eq:potential}
|g(\rho)|\le CG(\rho),
\een
if $g\in C^2$ with $g(\overline{\rho})=g'(\overline{\rho})=0$.

{\bf Step 2.}\,$H^1$ energy estimate

Multiply (\ref{eq:cNS}) by $\dot u$ and integrate the resulting equation over $\mathbf{R}^3$ to obtain
\ben\label{eq:eH1-1}
\int_{\R^3}\rho|\dot u|^2dx= \int_{\R^3}(-\dot u \cdot \nabla P+\mu\Delta u \cdot \dot u+\lambda \nabla \dv u\cdot \dot u)dx.
\een
By the continuity equation, we have
\ben\label{eq:renorm}
\pa_tP+\dv(uP)=\dv u(P-P'\rho).
\een
Then integration by parts yields that
\ben\label{eq:eH1-2}
&&\int_{\R^3}-\dot u \cdot \nabla P dx=\int_{\R^3}(\dv u_t(P-P(\overline{\rho}))-(u\cdot\nabla u)\cdot \nabla P)dx\nonumber\\
&&=\p_t\int_{\R^3}\dv u (P-P(\overline{\rho}))dx+\int_{\R^3}\big((P'\rho-P)(\dv u)^2+P\pa_iu^j\pa_ju^i\big)dx\nonumber\\
&&\leq \p_t\int\dv u(P-P(\overline{\rho})dx+C\|\nabla u\|^2_{2}.
\een
Integration by parts again gives
\ben
\mu\int_{\R^3}\Delta u \cdot \dot udx&=&-\f \mu2\p_t\int_{\R^3}|\na u|^2dx
-\mu\int_{\R^3}\pa_iu^j\pa_i(u^k\pa_ku^j)dx\nonumber\\
&\leq& -\f \mu2\p_t\int_{\R^3}|\na u|^2dx +C\int |\nabla u|^3 dx,
\een
and similarly,
\ben\label{eq:eh1-4}
\lambda \int_{\R^3}\nabla \dv u\cdot \dot u dx&=&-\f \lambda 2\p_t\int_{\R^3}|\dv u|^2dx
-\lambda\int_{\R^3}\pa_iu^j\pa_i(u^k\pa_ku^j)dx\nonumber\\
&\leq&-\f \lambda 2\p_t\int_{\R^3}|\dv u|^2dx+C\int_{\R^3}|\nabla u|^3dx.
\een
And by (\ref{eq:L2}) and (\ref{eq:potential}), we have
\beno
\int_{\R^3}\dv u(P-P(\overline{\rho})dx\leq \|\dv u\|_{2} \|P-P(\overline{\rho})\|_{2}
\leq CC_0^{1/2}\|\dv u\|_{2}.
\eeno
Plugging (\ref{eq:eH1-2})-(\ref{eq:eh1-4}) into (\ref{eq:eH1-1}) yields that
\ben\label{eq:A1}
A_1(T_0)\leq CC_0+C\int_0^{T_0}\int|\nabla u|^3 dxdt.
\een

{\bf Step 3.}\, $H^2$ energy estimate

We take the material derivative to the second equation of (\ref{eq:cNS}) to obtain
\beno
&&\rho \dot{u}_t+\rho u\cdot \nabla \dot{u}+\nabla P_{t}+\dv(\nabla P\otimes u)\\
&&\quad=\mu\big[\Delta u_t+\dv(\Delta u\otimes u)\big]+
\lambda\big[\nabla\dv u_{t}+\dv((\nabla\dv u)\otimes u)\big].
\eeno
Let us introduce some notations. We denote $\dv(f\otimes u)=\sum_{j=1}^3\partial_j(fu_j)$.
For two matrices $A=(a_{ij})_{3\times 3}$ and $B=(b_{ij})_{3\times 3}$, we use the notation $A:B=\sum_{i,j=1}^3a_{ij}b_{ij}$ and $AB$ is as usual the multiplication of matrix.

Multiply the above equation by $\sigma(t)\dot u$ and integrate by parts to get
\beno
&&\Big(\f{\sigma}{2}\int_{\R^3}\rho |\dot u|^2dx\Big)_t-\f12 \sigma'\int_{\R^3}\rho|\dot u|^2dx\\
&&=-\int_{\R^3}\sigma\dot u^j[\pa_j P_t+\dv (\pa_jPu)]dx
+\mu\sigma\int_{\R^3}\dot{u}\cdot\big(\Delta u_t+\dv(\Delta u\otimes u)\big)dx\\
&&\qquad+\lambda\sigma\int_{\R^3}\dot{u}\cdot\Big[\nabla\dv u_{t}+\dv((\nabla\dv u)\otimes u)\Big]dx.
\eeno
Now integration by parts and use (\ref{eq:renorm}) to obtain
\beno
&&-\int_{\R^3}\sigma  \dot u^j[\pa_j P_t+\dv (\pa_jPu)]dx
=-\int_{\R^3}\sigma[(u\cdot\nabla P+P'\rho\dv u)\pa_j\dot u^j-\pa_k\dot u^j\pa_j P u^k]dx\\
&&=-\int_{\R^3}\sigma[P'\rho \dv u \pa_j \dot u^j-P\pa_k(u^k\pa_j \dot u^j)+P\pa_j(\pa_k\dot u^j u^k)]dx\\
&&\leq C\sigma\|\nabla u\|_{2}\|\nabla \dot u\|_{2}.
\eeno
Integrate by parts again to obtain
\beno
&&\mu\int_{\R^3}\dot{u}\cdot\big(\Delta
u_t+\dv(\Delta u\otimes u)\big)dx=-\mu\int_{\R^3}\left[\nabla\dot{u}:\nabla u_t+ u\otimes\Delta u:\nabla \dot{u}\right]dx\\
&&=-\mu\int_{\R^3}\Big[|\nabla\dot{u}|^2-\nabla(u\cdot\nabla u):\nabla\dot{u}+ u\otimes\Delta u:\nabla \dot{u}\Big]dx\\
&&=-\mu\int_{\R^3}\Big[|\nabla\dot{u}|^2-\big((\nabla u\nabla u)+(u\cdot\nabla) \nabla u\big):\nabla\dot{u}-\nabla(u\cdot\nabla\dot{u}):\nabla u\Big]dx\\
&&=-\mu\int_{\R^3}\Big[|\nabla\dot{u}|^2-(\nabla u\nabla u):\nabla\dot{u}-\dv(\nabla u\otimes u):\nabla\dot{u}-(\nabla u\nabla\dot{u}):\nabla u-((u\cdot\nabla)\nabla\dot{u}):\nabla u\Big]dx\\
&&=-\mu\int_{\R^3}\Big[|\nabla\dot{u}|^2-(\nabla u\nabla u):\nabla\dot{u}+((u\cdot\nabla)\nabla\dot{u}):\nabla u-(\nabla u\nabla\dot{u}):\nabla u-((u\cdot\nabla)\nabla\dot{u}):\nabla u\Big]dx\\
&&\leq \int_{\R^3}\Big[-\f{3\mu}{4}|\nabla\dot{u}|^2+C|\nabla u|^4\Big]dx.
\eeno
To estimate the last term, note that
\beno
&&\dv((\nabla\dv u)\otimes u)=\nabla(u\cdot\nabla\dv u)-\dv(\dv u\nabla\otimes u)+\nabla(\dv u)^2,\\
&&\dv\dot{u}=\dv u_t+\dv(u\cdot\nabla u)=\dv u_t+u\cdot\nabla\dv u+\nabla u:(\nabla u)'.
\eeno
Here $A'$ means the transpose of matrix $A$. We have
\beno
&&\lambda\int_{\R^3}\dot{u}\cdot\Big[\nabla\dv u_{t}+\dv((\nabla\dv u)\otimes u)\Big]dx\\
&&=-\lambda\int_{\R^3}\Big[\dv\dot{u}\dv u_t+\dv\dot{u}(u\cdot\nabla\dv u)-\dv u(\nabla\dot{u})':\nabla u+\dv\dot{u}(\dv u)^2\Big]dx\\
&&=-\lambda\int_{\R^3}\Big[|\dv\dot{u}|^2-\dv\dot{u}\nabla u:(\nabla u)'-\dv u(\nabla\dot{u})':\nabla u+\dv\dot{u}(\dv u)^2\Big]dx\\
&&\leq\int_{\R^3}\Big[-\f{\lambda}{2}|\dv\dot{u}|^2+\f{1}{4} |\nabla\dot{u}|^2+C|\nabla u|^4\Big]dx.
\eeno
Summing up the above estimates, we obtain
\ben\label{eq:A2}
A_2(T_0)\leq CC_0+CA_1(T_0)+C\int^{T_0}_0\int_{\R^3}\sigma|\nabla u|^4dxdt,
\een
which along with (\ref{eq:A1}) gives
\ben\label{eq:A}
A_1(T_0)+A_2(T_0)\leq CC_0+C\int^{T_0}_0\int_{\R^3}\sigma|\nabla u|^4dxdt+C\int^{T_0}_0\int_{\R^3}|\nabla u|^3dxdt.
\een

It remains to bound the right hand side. For this purpose, we set
$$
F=(\mu+\lambda)\dv u-P(\rho)+P(\overline{\rho}),\quad \omega=\nabla\times u.
$$
By (\ref{eq:renorm}), we have
\beno
&&\pa_t(P-P(\overline{\rho}))+u\cdot\nabla (P-P(\overline{\rho}))\\
&&=\f{P'\rho}{P}\big[\f{1}{\mu+\lambda}F(P(\rho)-P(\overline{\rho}))
+\f{1}{\mu+\lambda}(P-P(\overline{\rho}))^2+P(\overline{\rho})\dv u\big],
\eeno
which implies that
\beno
&&C\|P(\rho)-P(\overline{\rho})\|^4_{4}\\
&&\leq \p_t\int(P-P(\overline{\rho}))^3dx+\delta\|P(\rho)-P(\overline{\rho})\|^4_{4}+
C(\|F\|^4_{4}+\|\nabla u\|^2_{2}).
\eeno
Taking $\delta$ small and integrating the above equation on $[0,T_0]$, we get
\beno
\int^{T_0}_0\sigma\|P(\rho)-P(\overline{\rho})\|^4_{L^4}dt
&\leq& C\sup_{t\in[0,T]}\|P(\rho)-P(\overline{\rho})\|^3_{L^3}
+C\int^{\sigma(T)}_0\|P(\rho)-P(\overline{\rho})\|^3_{L^3}dt\\
&&+C\int^{T_0}_0\sigma\|F\|^4_{4}dt+CC_0\\
&\leq& CC_0 +C\int^{T_0}_0 \sigma\|F\|^4_{L^4}dt,
\eeno
from which and Lemma \ref{lem:1}, we infer that
\beno
&&\int^{T_0}_0\int_{\R^3}\sigma |\nabla u|^4dxdt\\
&&\leq  C\int^{T_0}_0 \sigma (\|F\|^4_{4}+\|\omega\|^4_{4})dt
+C\int^{T_0}_0 \sigma\|P(\rho)-P(\overline{\rho})\|^4_{4}dt\\
&&\leq C\int^{T_0}_0 \sigma(\|\nabla u\|_{2}+\|P(\rho)-P(\overline{\rho})\|_{2})\|\rho \dot u\|^3_{2}dt
+ CC_0\\
&&\leq C\sup_{t\in [0,T_0]}(\sigma^{1/2}\|\rho^{1/2}\dot u\|_{2}(\|\nabla u\|_{2}
+C_0^{1/2}))\int^{T_0}_0 \|\rho^{1/2} \dot u\|^2_{2}dt
+ CC_0\\
&&\leq C(A_1^{1/2}+C_0^{1/2})A_2^{1/2}A_1(T)+C_0
\leq CC_0,
\eeno
and
\beno
&&\int^{T_0}_0\int_{\R^3}|\nabla u|^3dxdt\\
&&\leq \int^{\sigma(T_0)}_0 \int_{\R^3}|\nabla u|^3dxdt+\int^{T_0}_{\sigma(T_0)}\int_{\R^3}|\nabla u|^3dxdt\\
&&\leq \int^{\sigma(T_0)}_0 \|\nabla u\|_{2}^{\f3 2}(\|\rho\dot u\|^{\f32}_{2}+\|P-P(\overline{\rho})\|^{\f32}_{2})dxdt
+\int^{T_0}_{\sigma(T_0)}\int_{\R^3}\big(|\nabla u|^4+|\nabla u|^2\big)dxdt\\
&&\leq C\sup_{t\in[0,\sigma(T_0)]} \|\nabla u\|_{2}^{\f32} \int^{\sigma(T_0)}_0 \|\rho\dot u\|^{\f32}_{2}
+CC_0\\
&&\leq CA_1^{\f32}+CC_0\leq  CC_0^{\f34}.
\eeno
Plugging them into (\ref{eq:A}) and taking $\e$ small enough
depending on $\mu,\lambda,\overline{\rho},c_0$, we conclude
\ben\label{eq:A1+A2}
A_1(T_0)+A_2(T_0)\leq \f 12C_0^{\f12}.
\een

{\bf Step 4.}\,\,Lower and upper bound of the density

Set $L=\log(\rho)$, which  satisfies
$$
(2\mu+\lambda)\dot{L}+(P(\rho)-P(\overline{\rho}))=-F.
$$
For $0<t<\sigma(T_0)$, we have
\begin{eqnarray*}
\int_0^{\sigma(T_0)}\|F\|_{L^\infty}d{t}
&\le&C\int_0^{\sigma(T_0)}\|F\|_{L^6}^{\frac{1}{2}}\|\nabla{F}\|_{L^6}^{\f12}d{t}\\
&\le&C\int_0^{\sigma(T_0)}\|\rho\dot{u}\|_{L^2}^{\frac{1}{2}}\|\rho\dot{u}\|_{L^6}^{\f12}d{t}\\
&\le&C\int_0^{\sigma(T_0)}\|\rho^{1/2}\dot{u}\|_{L^2}^{\frac{1}{2}}\|\nabla{\dot{u}}\|_{L^2}^{\f12}d{t}\\
&=&C\int_0^{\sigma(T_0)}(t^{-1/2})^{\f12}(\|\rho^{1/2}\dot{u}\|_{L^2}^2)^{\f14}
(t\|\nabla{\dot{u}}\|_{L^2}^2)^{\f14}d{t}\\
&\le& CC_0^{\f 14},
\end{eqnarray*}
which implies that for $t\le\sigma(T_0)$,
\beno
\inf(\log \rho_0(x))-CC_0^{\f 14}-Ct\le \log\rho(t,x)\le \sup(\log\rho_0(x))+CC_0^\f14+Ct,
\eeno
hence, one can choose $\varepsilon_0,\tau$ sufficiently small such that for $t\le\tau\le \sigma(T_0)$,
\ben\label{eq:density-point}
\f3 4 c_0<\rho(t,x)<\f 32c_0^{-1}.
\een
For $\tau\le t_1\le t_2\le T_0$, we have
\begin{eqnarray*}
\int_{t_1}^{t_2}\|F\|_{L^\infty}d{t}
&\le& C\int_{t_1}^{t_2}\|\rho^{1/2}\dot{u}\|_{L^2}^{\frac{1}{2}}\|\nabla{\dot{u}}\|_{L^2}^{\f12}d{t}\\
&\le& C(t_2-t_1)^{1/2}\int_{t_1}^{t_2}(\|\rho^{1/2}\dot{u}\|_{L^2}^2+
\|\nabla{\dot{u}}\|_{L^2}^2)d{t}\\
&\le&\eta\big(t_2-t_1\big)+CC_0
\end{eqnarray*}
for any $\eta>0$. Now let $T^*$ be a maximal time so that (\ref{eq:density-point}) holds.
If $T^*\ge T$, then we are done. Otherwise, we have
\beno
\rho(T^*,x(T^*))=\f 32c_0^{-1}\quad \textrm{or} \quad \rho(T^*,x(T^*))=\f 34c_0,
\eeno
where the curve $x(t)$ is defined by
\beno
\f d {dt}x(t)=u(t,x(t)),\quad x(0)=x.
\eeno
If $\rho(T^*,x(T^*))=\f 32c_0^{-1}$, then there exits $T^*>T_1\ge \tau$ such that
\beno
\rho(T_1,x(T_1))=c_0^{-1},\quad c_0^{-1}\le \rho(t,x(t))\le \f 32 c_0^{-1}\quad \textrm{for}\quad t\in [T_1,T^*].
\eeno
Then we conclude by (\ref{ass:pressure}) that
\beno
L(T^*,x(T^*))-L(T_1,x(T_1))\le CC_0,
\eeno
if we choose $\eta$ small enough. If we choose $\varepsilon_0$ sufficiently small, this in turn
implies that
\beno
\rho(T^*,x(T^*))<\f 32c_0^{-1},
\eeno
which contradicts with the definition of $T^*$. Similar argument yields that
\beno
\rho(T^*,x(T^*))>\f 34c_0,
\eeno
if $\rho(T^*,x(T^*))=\f 34c_0$. Hence, we have
\ben\label{eq:density-f}
\f3 4 c_0<\rho(t,x)<\f 32c_0^{-1},\quad \textrm{for}\quad 0\le t\le T_0.
\een

With (\ref{eq:A1+A2}) and (\ref{eq:density-f}),
we conclude the proof by continuity argument.\ef

\begin{Lemma}\label{lem:1}
For any $p\in[2,6]$, there exists a constant $C$ depending only on $\mu,\lambda$ such that
\beno
&&\|\nabla u\|_{p}\leq C\big(\|F\|_{p}+\|\omega\|_{p}+\|P(\rho)-P(\overline{\rho})\|_{p}\big),\\
&&\|\nabla u\|_{p}\leq C\|\nabla u\|_{2}^{\f{6-p}{2p}}
\big(\|\rho\dot u\|_{2}+\|P(\rho)-P(\overline{\rho})\|_{6}\big)^{\f {3p-6}{2p}}.
\eeno
\end{Lemma}

\no{\bf Proof.}\, By the elliptic estimate and the definition of $F$, we have
\beno
\|\nabla u\|_{p}\leq C\big(\|\dv u\|_{p}+\|\omega\|_{p}\big)
\leq C\big(\|F\|_{p}+\|\omega\|_{p}+\|P(\rho)-P(\overline{\rho})\|_{p}\big).
\eeno
Noting that
$$\Delta F=\dv(\rho \dot u),\quad\mu\Delta\omega=\nabla\times(\rho \dot u),$$
we have by the elliptic estimate that
\beno
\|\nabla F\|_{2}+\|\nabla \omega\|_{2}\leq C\|\rho \dot u\|_{2}.
\eeno
Then by H\"{o}lder inequality and Sobolev inequality, we get
\beno
\|\nabla u\|_{p}&\leq& C\|\nabla u\|^{\f {6-p}{2p}}_{2}\|\nabla u\|^{\f {3p-6}{2p}}_{6}\\
&\leq& C\|\nabla u\|^{\f {6-p}{2p}}_{2}
\big(\|F\|_{6}+\|\omega\|_{6}+\|P(\rho)-P(\overline{\rho})\|_{6})^{\f{3p-6}{2p}}\\
&\leq& C\|\nabla u\|^{\f {6-p} {2p}}_{2}\big(\|\rho \dot u\|_{2}+\|P(\rho)-P(\overline{\rho})\|_{6})^{\f{3p-6}{2p}}.
\eeno
This finishes the proof of lemma.\ef

\setcounter{equation}{0}
\section{Continuation criterion}

In order to extend a local solution to a global one,
we need to establish a continuation criterion of smooth solution.
The proof is motivated by \cite{Sun}.

\begin{Theorem}\label{thm:con}
Let $(\rho,u)$ be a solution of (\ref{eq:cNS}) satisfying
\beno
\rho(0)>0,\quad \rho-\overline{\rho} \in C([0,T];H^2),\quad u\in C([0,T];H^2)\cap L^2(0,T;H^3).
\eeno
Let $T^*$ be the maximal existence time of the solution. If $T^*<+\infty$, then it is necessary
\ben
\limsup_{t\uparrow T^*}\big(\|\rho(t)\|_{\infty}+\|u(t)\|_q\big)=+\infty
\een
for any $q>3$.
\end{Theorem}

\begin{Remark}
Compared with Theorem 1.3 in \cite{Sun}, the restriction on the viscosity coefficient $(\lambda,\mu)$ is removed,
but the price to pay is to impose the condition on the velocity $u$. Due to $\rho(0)>0$,
the compatibility condition on the initial data is also removed.
\end{Remark}

\no{\bf Proof.}\,We use the contradiction argument. Hence, assume that $T^*<\infty$ and
\ben\label{eq:condition}
\sup\limits_{t\in [0,T^*)}\big(\|\rho(t)\|_{\infty}+\|u(t)\|_q\big)=M<+\infty.
\een
In what follows, we denote $C$
by a constant depending on $T, M, \|u_0\|_{H^2}$ and $\|\rho_0-\overline{\rho}\|_{H^2}$.

First of all, we have
\beno
\int_{\R^3}\f 12\rho(t,x)|u(t,x)|^2+G(\rho(t,x))dx+\int_0^t\int_{\R^3}|\na u|^2dxdt\le C.
\eeno
This gives by (\ref{eq:condition}) and (\ref{eq:potential}) that for any $r\in [2,\infty]$,
\ben\label{eq:L2-energy}
\|\rho-\overline{\rho}\|_{L^\infty(0,T;L^r)}+\|\sqrt{\rho}u\|_{L^\infty(0,T;L^2)}
+\|\na u\|_{L^2(0,T;L^2)}\le C.
\een

Let $v=L^{-1}\na P(\rho)$ be a solution of the following elliptic system
\beno
Lv\eqdefa \mu\Delta v+\lambda\nabla \textrm{div} v=\nabla P(\rho).
\eeno
By elliptic estimate and (\ref{eq:L2-energy}), for $r\in [2,\infty)$
\ben\label{eq:v-Lr}
\|\na v\|_r\le C\|P(\rho)-P(\overline{\rho})\|_r\le C,
\quad\|\na^2 v\|_r\le C\|\na \rho\|_r.
\een

Introduce a new unknown $w=u-v$. Note that $\textrm{div} w=(\lambda+\mu)\textrm{div} u-P(\rho)+P(\overline{\rho})$ is a so called
the effective viscous flux introduced by Hoff \cite{Hoff-JDE}. It is easy to find that $w$ satisfies
\beno
\rho\partial_t w-\mu\Delta w-\lambda\nabla \textrm{div} w=\rho F,
\eeno
with
\beno
F=-u\cdot\nabla u+L^{-1}\nabla \textrm{div}\big(p(\rho)u\big)+L^{-1}\nabla\big((\rho
p'(\rho)-p(\rho))\textrm{div} u\big).
\eeno
Multiplying the above equation by $\partial_t w$ and integrating by parts, we obtain
\ben\label{eq:w-energy}
\frac{d}{dt}\int_{\R^3}\mu|\nabla
w|^2+\lambda|\textrm{div} w|^2dx+\f{1}{2}\int_{\R^3}\rho|\partial_t w|^2dx
\leq \f{1}{2}\|\sqrt{\rho}F\|^2_{2}.
\een
By Sobolev inequality, interpolation inequality, (\ref{eq:condition}) and (\ref{eq:v-Lr}), we get
\beno
\|\sqrt{\rho} u\cdot\nabla
u\|_{L^2(\Om)}
&\leq& C\|u\|_{q}\|\nabla u\|_{\f{2q}{q-2}}
\leq C\big(\|\nabla w\|_{\f{2q}{q-2}}+\|\nabla v\|_{\f{2q}{q-2}}\big)\\
&\leq& C\|\nabla w\|_{2}+\epsilon\|\nabla^2
w\|_{2}+C\\
&\le& C+C\|\na u\|_2+\epsilon\|\nabla^2 w\|_{2}
\eeno
for any $\epsilon>0$, and by elliptic estimate, (\ref{eq:L2-energy}) and (\ref{eq:v-Lr}),
\beno
\|\sqrt{\rho}
L^{-1}\nabla \textrm{div}[p(\rho)u]\|_{2}\leq
C\|p(\rho)u\|_{2}\leq C\|\sqrt{\rho}u\|_{2}\leq C,
\eeno
and by Sobolev inequality,
\beno
&&\|\sqrt{\rho} L^{-1}\nabla (\rho p'-p)\textrm{div} u\|_{2}\\
&&\le \|\rho-\overline{\rho}\|_3\| L^{-1}\nabla (\rho p'-p)\textrm{div} u\|_6
+C\|(\rho-\overline{\rho})\textrm{div} u\|_{\f 65}+C{\overline{\rho}}\|u\|_2\\
&&\le C\|\na u\|_2+C\|\sqrt{\rho}u\|_2\le C+C\|\na u\|_2,
\eeno
here we used
\beno
\overline{\rho}^2\int_{\R^3}|u|^2dx\le 2\int_{\R^3}|\rho-\overline{\rho}|^2|u|^2dx+2\int_{\R^3}\rho^2|u|^2dx.
\eeno
Hence,
\beno
\|\sqrt{\rho}F\|_{2}\leq C(1+\|\nabla u\|_{2})+\epsilon\|\nabla^2w\|_{2}.
\eeno
Noting that $Lw=\rho\partial_t w-\rho F$, we get by elliptic estimate that
\beno
\|\nabla^2w\|_{2}\leq C\big(\|\rho\partial_t
w\|_{2}+\|\rho F\|_{2}\big)\leq
C\big(\|\sqrt{\rho}\partial_t w\|_{2}+\|\sqrt{\rho}
F\|_{2}\big),
\eeno
which implies by taking $\epsilon$ small enough that
\beno
\|\sqrt{\rho}F\|_{2}\leq C(1+\|\nabla u\|_{2})+\f{1}{2}\|\sqrt{\rho}\partial_t
w\|_{2}.
\eeno
This along with (\ref{eq:w-energy}) gives
\beno
\|\na w\|_{L^\infty(0,T;L^2)}+\|\sqrt{\rho}\p_tw\|_{L^2(0,T;L^2)}+\|\na^2 w\|_{L^2(0,T;L^2)}\le C,
\eeno
which in turn implies by (\ref{eq:v-Lr}) and Sobolev inequality that
\ben\label{eq:u-h1}
\|\na u\|_{L^\infty(0,T;L^2)}+\|\na u\|_{_{L^2(0,T;L^r)}}\le C\quad \textrm{for} \quad r\in [2,6].
\een

Now let us turn to the high order energy estimate. Thanks to the proof of Step 3 in Theorem \ref{thm:hoff}, we have
\ben\label{eq:u-H2}
&&\f{d}{dt}\int_{\R^3}t\rho|\dot{u}|^2dx+\mu t\int_{\R^3}|\nabla\dot{u}|^2dx
+\lambda t\int_{\R^3}|\textrm{div}\dot{u}|^2dx\nonumber\\
&&\leq C+Ct\int_{\R^3}|\nabla{u}|^4dx+C\int_{\R^3}|\nabla u|^3dx.
\een
Noting that
\beno\label{eq:w}
\mu\Delta w+\lambda\nabla \textrm{div}w=\rho \dot{u},
\eeno
elliptic estimate yields that
\beno
\|\nabla^2 w\|_{2}\leq C\|\rho \dot{u}\|_{2}\leq C\|\sqrt{\rho}\dot{u}\|_{2}.
\eeno
Then by Sobolev inequality, (\ref{eq:v-Lr}) and (\ref{eq:u-h1}), we infer that
\beno
\|\nabla u\|^4_{4}&\leq&\|\nabla u\|_{2}\|\nabla
u\|^3_{6}\leq C\|\nabla u||^2_{6}\big(\|\nabla
w\|_{6}+\|\nabla v\|_{6}\big)\\
&\leq& C\|\nabla u||^2_{6}\big(1+\|\nabla^2
w\|_{2}\big)\leq C\|\nabla u\|^2_{6}\big(1+\|\sqrt{\rho}\dot{u}\|_{2}\big),\\
\|\na u\|_3^3&\leq& \|\na u\|_2^\f32\|\na u\|_6^\f32\le C\|\na u\|_6^\f32.
\eeno
Plugging them into (\ref{eq:u-H2}) and noting $||\nabla u(t)||_{6}\in L^2(0,T)$ by (\ref{eq:u-h1}),
we deduce by Gronwall's inequality that
\ben\label{eq:du-h2}
t\int_{\R^3}\rho|\dot{u}|^2dx+\int^T_0\int_{\R^3}t|\nabla \dot{u}|^2dxdt\leq C.
\een
from which, elliptic estimate and Sobolev inequality, it follows that
\ben\label{eq:w-Lr}
\|\nabla^2w\|_{L^2(0,T;L^r)}\leq C\|\rho\dot u\|_{L^2(0,T;L^r)}\le C\quad \textrm{for} \quad r\in [2,6].
\een

To obtain the high order estimate for $u$, we need to use the continuity equation.
Take the derivative with respect to $x$ for the first equation of (\ref{eq:cNS}) to obtain
\beno
\partial_t \nabla\rho+(u\cdot\nabla)\nabla \rho+\nabla u\nabla\rho+\textrm{div} u\nabla\rho+\rho\nabla \textrm{div} u=0.
\eeno
For $r\in (3,6]$, making $L^r$ energy estimate yields that
\beno
&&\frac{d}{dt}\int_{\R^3}|\nabla\rho|^rdx\\
&&\leq C\int_{\R^3}|\nabla u||\nabla\rho|^rdx
+C\int_{\R^3}\rho|\nabla \textrm{div} u||\nabla\rho|^{r-1}dx\\
&&\leq C\|\nabla
u\|_{\infty}\|\nabla\rho\|^r_{r}+C\|\nabla^2u\|_{r}\|\nabla\rho\|^{r-1}_{r}\\
&&\leq C\big(\|\nabla w\|_{\infty}+\|\nabla
v\|_{\infty}\big)\|\nabla\rho\|^r_{r}+C\big(\|\nabla^2w\|_{r}
+\|\nabla^2v\|_{r}\big)\|\nabla\rho\|^{r-1}_{r}.
\eeno
Then by (\ref{eq:v-Lr}) and Sobolev inequality, we get
\beno
\frac{d}{dt}\int_{\R^3}|\nabla\rho|^rdx
\le C\big(1+\|\nabla v\|_{\infty}+\|\nabla^2 w\|_{r}\big)\|\nabla\rho\|^r_r,
\eeno
and by Lemma \ref{lem:log} and elliptic estimate,
\beno
\|\nabla v\|_{\infty}&\leq&
C\big(\|\na v\|_2+\|\nabla v\|_{\dot B^{0}_{\infty,\infty}}\ln(e+\|\nabla v\|_{W^{1,r}})\big)\\
&\leq& C\big(\|\rho-\overline{\rho}\|_2+\|\rho\|_{\infty}\ln(e+\|\nabla\rho\|_{r})\big),
\eeno
from which and (\ref{eq:w-Lr}), we deduce by Gronwall's inequality that
\ben\label{eq:density-w1r}
\|\na\rho\|_r\le C,
\een
which along with (\ref{eq:w-Lr}) and (\ref{eq:v-Lr}) gives
\ben\label{eq:u-w2r}
\|\na^2 u\|_{L^2(0,T;L^r)}\le \|\na^2 w\|_{L^2(0,T;L^r)}+\|\na^2 v\|_{L^2(0,T;L^r)}\le C.
\een
This in turn implies (\ref{eq:density-w1r}) also holds for $r\in [2,6]$.

Now we are in position to give $H^2$ estimate of $(\rho,u)$. First of all, we have
\ben\label{eq:d-h2}
\f d {dt}\|\rho-\overline{\rho}\|_{H^2}^2\le C(1+\|\na u\|_{\infty})\|\rho-\overline{\rho}\|_{H^2}^2+C\|\na^3u\|_{2}^2.
\een
While by (\ref{eq:density-w1r}),
\beno
\|\na (\rho\dot u)\|_2
\le\|\na \rho\|_3\|\dot u\|_6+C\|\na \dot u\|_2\le C\|\na \dot u\|_2,
\eeno
this along with elliptic estimate and (\ref{eq:du-h2}) implies that
\beno
\|\na^2 w(t)\|_2^2+\int_0^t\|\na^3 w(\tau)\|_2^2d\tau\le C.
\eeno
Indeed, we have
\beno
\int_0^t\|\na^3 w(\tau)\|_2^2d\tau&\le& \int_0^\delta\|\na^3 w(\tau)\|_2^2d\tau+\int_\delta^t\|\na^3 w(\tau)\|_2^2d\tau\\
&\le& C+\delta^{-1}\int_0^t\tau\|\na^3 w(\tau)\|_2^2d\tau\le C.
\eeno
Hence by (\ref{eq:v-Lr}), we get
\ben\label{eq:u-H22}
\|\na^2 u(t)\|_2^2+\int_0^t\|\na^3 u(\tau)\|_2^2d\tau\le C+C\int_0^t\|\na^2\rho(\tau)\|_2^2d\tau.
\een
Summing up (\ref{eq:d-h2}) and (\ref{eq:u-H22}), we conclude by Gronwall's inequality that for $0\le t<T^*$,
\beno
\|u(t)\|_{H^2}+\|\rho(t)-\overline{\rho}\|_{H^2}\le C.
\eeno
This ensures that the solution can be continued after $t=T^*$.\ef

\begin{Lemma}\label{lem:log}
Let $r>3$ and $u\in W^{1,r}$. Then there holds
\beno
\|u\|_\infty\le C\|u\|_2+C\|u\|_{\dot B^0_{\infty,\infty}}\ln(e+\|u\|_{W^{1,r}}).
\eeno
\end{Lemma}

\no{\bf Proof.}\,The proof is standard. For the reader's convenience, we present a proof.
We decompose $u$ into
\beno
u=S_{-N}u+\sum_{j=-N}^{N}\Delta_ju+\sum_{j\ge N+1}\Delta_j u,
\eeno
from which and Lemma \ref{Lem:Bernstein}, we infer that
\beno
\|u\|_{\infty}&\le& 2^{-\f 32N}\|u\|_{2}+(2N+1)\|u\|_{\dot B^{0}_{\infty,\infty}}+C\sum_{j\ge N}2^{\f{3}rj}\|\Delta_j u\|_{r}\\
&\le& 2^{-\f 32N}\|u\|_{2}+(2N+1)\|u\|_{\dot B^{0}_{\infty,\infty}}+C2^{-(1-\f{3}r)N}\|u\|_{W^{1,r}}.
\eeno
Taking $N\in \N$ such that
\beno
2^{-(1-\f 3r)N}\|u\|_{W^{1,r}}\sim 1,
\eeno
the desired estimate follows easily.\ef

\setcounter{equation}{0}

\section{Proof of Theorem \ref{thm:global}}

Theorem \ref{thm:local}, Proposition \ref{prop:low} and Proposition \ref{prop:high} ensure that
there exists a solution $(\rho,u)$ of (\ref{eq:cNS}) satisfying
\beno
&&\f {c_0} 2\le \rho\le 2c_0^{-1},\quad \rho-\overline{\rho}\in \widetilde{L}^\infty(0,T;\dot B^{\f 3p}_{p,1}\cap\dot B^{s}_{2,2}),\\
&&u\in \widetilde{L}^\infty(0,T;\dot B^{\f 3p-1}_{p,1}\cap\dot B^{s-1}_{2,2})
\cap \widetilde{L}^1(0,T;\dot B^{\f 3p+1}_{p,1}\cap\dot B^{s+1}_{2,2}).
\eeno
Moreover, there hold
\ben
&&\|\rho-\overline{\rho}\|_{\widetilde{L}^\infty_T(\dot B^{1-\delta}_{2,2})}
\leq C\big(\|\rho_0-\overline{\rho}\|_{\dot H^{1-\delta}}+\|u_0\|_{\dot H^{-\delta}}\big),\label{eq:solution-d}\\
&&\|u\|_{\widetilde{L}^\infty_T(\dot B^{-\delta}_{2,2})}\le C(1+E_0)\bigl(\|u_0\|_{\dot
H^{-\delta}}+T\|\rho_0-\overline{\rho}\|_{\dot H^{1-\delta}}\bigr),\label{eq:solution-u1}\\
&&\|u\|_{\widetilde{L}^1_T(\dot B^{2-\delta}_{2,2})}+\|u\|_{\widetilde{L}^2_T(\dot B^{1-\delta}_{2,2})}
\le C\bigl(\|u_0\|_{\dot
H^{-\delta}}+T\|\rho_0-\overline{\rho}\|_{\dot H^{1-\delta}}\bigr),\label{eq:solution-u2}\\
&&\|u\|_{\widetilde{L}^1_T(\dot B^{ \frac{3}{p}+1}_{p,1})}+\|u\|_{\widetilde{L}^2_T(\dot B^{ \frac{3}{p}}_{p,1})}
\le \eta_2\le 1.\label{eq:solution-u3}
\een
and by Remark \ref{rem:interval}, the existence time $T$ has a lower bound
\beno
T\ge\f 1 {C\big(1+\|\rho_0-\overline{\rho}\|_{\dot B^{\f 3p}_{p,1}\cap H^s}\big)^{(s+\f12)/(s-\f32)}}=T_1.
\eeno
Here and what follows, the constant $C$ depends only on $\lambda,\mu, c_0, \overline{\rho}, p, s$.
Due to $\|u_0\|_{\dot H^{\delta}}\le c_2$, we get by (\ref{eq:solution-u1}) and (\ref{eq:solution-u2}) that
\beno
&&\|u\|_{\widetilde{L}^\infty_T(\dot B^{-\delta}_{2,2})}\le C(1+E_0)\big(c_2+T\|\rho_0-\overline{\rho}\|_{\dot H^{1-\delta}}\big),\\
&&\|u\|_{\widetilde{L}^1_T(\dot B^{2-\delta}_{2,2})}+\|u\|_{\widetilde{L}^2_T(\dot B^{1-\delta}_{2,2})}
\le C\big(c_2+T\|\rho_0-\overline{\rho}\|_{\dot H^{1-\delta}}\big).
\eeno
Then for $r\in (1,2)$, we get by Lemma \ref{Lem:Bernstein} that
\beno
&&\|\na u\|_{L^r_T(L^2)}\le \sum_{j\in \Z}2^j\|\Delta_j u\|_{L^r_T(L^2)}\\
&&\le \sum_{j\le 0}2^jT^\f 1r\|\Delta_j u\|_{L^\infty_T(L^2)}
+\sum_{j>0}2^j\|\Delta_j u\|_{L^\infty_T(L^2)}^{1-\f 1r}\|\Delta_j u\|_{L^1_T(L^2)}^{\f 1r}\\
&&\le C\big(\|u\|_{\widetilde{L}^\infty_T(\dot B^{-\delta}_{2,2})}
+\|u\|_{\widetilde{L}^1_T(\dot B^{2-\delta}_{2,2})}\big)\\
&&\le C(1+E_0)\big(c_2+T\|\rho_0-\overline{\rho}\|_{\dot H^{1-\delta}}\big),
\eeno
and similarly,
\beno
\|u\|_{L^r_T(L^2)}&\le& C\big(T^\f 1r\|u\|_{\widetilde{L}^\infty_T(\dot B^{-\delta}_{2,2})}
+\|\na u\|_{L^r_T(L^2)}\big)\\
&\le& C(1+E_0)\big(c_2+T\|\rho_0-\overline{\rho}\|_{\dot H^{1-\delta}}\big).
\eeno
Hence, there exits a point $t_0\in (0,T)$ such that
\ben\label{eq:u-small}
\|u(t_0)\|_{H^1}\le C(1+E_0)\big(c_2/T^\f23+T^{\f 13}\|\rho_0-\overline{\rho}\|_{\dot H^{1-\delta}}\big).
\een

Recall that the density $\rho$ satisfies
\beno
\p_t(\rho-\overline{\rho})+u\cdot\na(\rho-\overline{\rho})=-\textrm{div}u\rho.
\eeno
Making $L^2$ energy estimate, we get
\beno
\p_t\|\rho-\overline{\rho}\|_2\le \|\textrm{div} u\|_\infty\|\rho-\overline{\rho}\|_2+C\|\textrm{div}u\|_2.
\eeno
Then by Gronwall's inequality and (\ref{eq:solution-u3}), we get
\ben
\|\rho(t)-\overline{\rho}\|_2&\le& C\big(\|\rho_0-\overline{\rho}\|_2+T^\f12\|\na u\|_{L^2_T(L^2)}\big)\nonumber\\
&\le& C\big(c_1+T^\f12\|u_0\|_2\big).\label{eq:d-small}
\een
Here we used
\beno
\|\na u\|_{L^2_T(L^2)}\leq C\big(\|\rho_0-\overline{\rho}\|_{2}+
\|u_0\|_{2}\big),
\eeno
which follows from the energy inequality (\ref{eq:L2}). Set
\beno
T_2=\f {\varepsilon_0^2} {C\big(1+\|\rho_0-\overline{\rho}\|_{\dot B^{\f 3p}_{p,1}\cap H^s}+\|u_0\|_2\big)^{6}},
\eeno
and $T_*=\min(T_1,T_2)$. We choose $c_1, c_2$ small enough so that
\beno
Cc_1\le \f {\varepsilon_0} 4,\quad C(1+E_0)c_2/T^\f23_*\le \f {\varepsilon_0} 4,
\eeno
i.e, $c_2\le \varepsilon_0T^\f23_*/(4C(1+E_0))$.
Then it follows from (\ref{eq:u-small}) and (\ref{eq:d-small}) that
\beno
\|a(t_0)\|_2+\|u(t_0)\|_{H^1}\le \varepsilon_0.
\eeno
Theorem \ref{thm:hoff} ensures that
\beno
\f {c_0} 4\le\rho\le 4c_0^{-1},\quad \|u(t)\|_{L^\infty_T(L^q)}\le C.
\eeno
So, the solution can be extended to a global one by Theorem \ref{thm:con}.\ef

\noindent {\bf Acknowledgments.}
This work is carried out while the third author is a
long term visitor at Beijing international Mathematical center of  research(BIMCR).
The hospitality and support of BIMCR are graciously acknowledged.
Zhifei Zhang was partially supported by the NSF of China under grants 10990013 and 11071007.

\end{document}